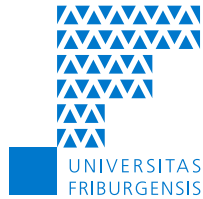
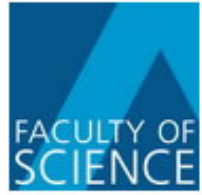

# On the Conjugacy Problem in Groups and its Variants

A thesis submitted for the degree of
Master of Science in Mathematics
University of Fribourg

Michèle Feltz
Department of Mathematics
University of Fribourg
Switzerland

Supervisor: Dr. Laura Ciobanu

February 2010

**Abstract**

This thesis deals with the conjugacy problem in groups and its twisted variants. We analyze recent results by Bogopolski, Martino, Maslakova and Ventura on the twisted conjugacy problem in free groups and its implication for the conjugacy problem in free-by-cyclic groups and some further group extensions. We also consider the doubly-twisted conjugacy problem in free groups. Staecker has developed an algorithm for deciding doubly-twisted conjugacy relations in the case where the involved homomorphisms satisfy a certain remnant inequality. We show how a similar condition affects the equalizer subgroup and raise new questions regarding this subgroup. As an application we discuss the Shpilrain-Ushakov authentication scheme based on the doubly-twisted conjugacy search problem in matrix semigroups over truncated polynomials over finite fields.

Part of this thesis is devoted to the implementation and testing of some of the previously mentioned concepts in the GAP programming language.

**Acknowledgements**

I would like to thank my supervisor Dr. Laura Ciobanu for suggesting a highly interesting topic in a challenging field, for many helpful discussions and for constructive feedbacks on drafts of this work. My interest in group theoretical decision problems and their potential application to cryptography arose from a seminar on the word problem that I worked on under the guidance of Dr. Laura Ciobanu. Her commitment and her enthusiastic approach to research have always been a source of motivation to me.

I owe many thanks to my family and friends for their support during my studies in Fribourg.

# Contents

# Chapter 1

# Introduction

This thesis is intended to investigate recent progress in the field of group theoretical decision problems. These are problems of the following nature: given a property $P$ and an object $O$, find out if it is algorithmically decidable whether the object $O$ has the property $P$ or not.

In 1912, Max Dehn raised the following three fundamental decision problems that lie at the base of combinatorial group theory (see [26] for a survey).

- The *Word Problem* for a finitely presented group $G$ is the algorithmic problem of deciding, for an arbitrary word $\omega \in G$, whether or not $\omega =_G 1$.

- The *Conjugacy Problem* for a finitely presented group $G$ refers to the algorithmic problem of determining, for two arbitrary words $u, v \in G$, whether or not $u$ and $v$ are conjugate in $G$, that is, whether or not there exists an element $x \in G$ such that $x^{-1}ux =_G v$.

- The *Isomorphism Problem* is the decision problem of determining whether two finite group presentations present isomorphic groups.

Of special interest to us will be the conjugacy problem and its generalizations, namely the *twisted conjugacy problem* and the *doubly-twisted conjugacy problem*. These variants involve homomorphisms acting on the conjugating element.

All of the previously mentioned decision problems are unsolvable for general group presentations. So it is an interesting question to ask whether they might be solvable for classes of finitely generated groups enjoying some particular algebraic property (for instance, being free, abelian or hyperbolic). It has turned out that solving the twisted conjugacy problem in simpler groups is essential for studying the standard conjugacy problem in more complicated groups such as free-by-cyclic or free-by-free groups. Our main concern in the first part of this thesis is to analyze these developments. The twisted conjugacy problem will then naturally lead to the doubly-twisted conjugacy problem, which is a generalization of the former. Techniques and algorithms have been developed to decide doubly-twisted conjugacy relations in free groups by imposing *remnant* conditions on the involved homomorphisms. In particular, we investigate how these techniques affect the equalizer subgroup. Moreover, we briefly address the issue of whether there exists a relation between deciding doubly-twisted conjugacy relations and membership in the equalizer subgroup.

It is worth mentioning that, just as for Dehn's fundamental problems, motivation for studying the twisted variants of the conjugacy problem can be found in algebraic topology. Twisted conjugacy plays an important role in Nielsen fixed point theory (see [19]). This theory is concerned with the study of the *Nielsen number* of a selfmap $f : X \to X$ on a topological space $X$. The Nielsen number of $f$ is an algebraic invariant that gives a lower bound for the minimum number of fixed points in the homotopy class of $f$. Wagner showed in [36] that there exists an algorithm to compute the Nielsen number of a selfmap $f$ on a hyperbolic surface with boundary if the induced endomorphism $f_\#$ on the fundamental group, which is a finitely generated free group, has remnant (by considering the twisted conjugacy problem in the fundamental group). Nielsen coincidence theory (see [13]) focuses on the coincidence set $Coin(f, g) = \{x \in X \mid f(x) = g(x)\}$ of a pair of mappings $f, g : X \to Y$ and studies the *Nielsen coincidence number*, a lower bound for the minimal number of coincidence points considered in the homotopy classes of $f$ and $g$. For compact surfaces with boundary $X, Y$, deciding whether two coincidence points are in the same *Nielsen class* (meaning, that they can be combined by a homotopy) can be reduced to solving a doubly-twisted conjugacy problem in the fundamental groups, which are finitely generated free, by using the induced homomorphisms $f_\#, g_\#$.

In recent years there has been considerable interest in public-key cryptographic schemes based on the hardness of the conjugacy search problem in specific non-commutative groups (see [27] for a survey). The conjugacy search problem in a group $G$ is the problem of recovering a particular element $x \in G$ from two given elements $g \in G$ and $h = x^{-1}gx$. It can be considered as a generalization of the *discrete logarithm problem*. Several schemes using braid groups as platform groups have been proposed. However, it seems that the conjugacy search problem in braid groups is unlikely to provide a sufficient level of security due to heuristic attacks, in particular *length-based attacks*. In this work we will discuss a group-based authentication scheme which is based on the doubly-twisted conjugacy search problem in a semigroup of matrices over a finite commutative ring and we will see that this scheme is not secure against some heuristic algebraic attack. Thus, it remains a challenge to find hard problems in combinatorial group theory, as well as appropriate platform groups, that will allow us to develop provably secure cryptographic schemes.

An outline of this thesis follows.

In Chapter 2 we deal with the solvability of the twisted conjugacy problem in free groups and its consequence for free-by-cyclic groups. We point out how these techniques can be extended to larger families of groups. In particular, we explain how to construct [free abelian]-by-free groups with unsolvable conjugacy problem.

In Chapter 3 we explore remnant conditions on homomorphisms that have an impact on deciding doubly-twisted conjugacy relations. Random elements that belong to the same doubly-twisted conjugacy class have been generated using the GAP programming language. Moreover, we deduce a new condition forcing the equalizer subgroup to be trivial. Also some new open questions that might be of interest for future research are discussed.

In Chapter 4 we introduce a public-key authentication scheme based on the doubly-twisted conjugacy search problem in a matrix semigroup over truncated one-variable polynomials over the finite field of two elements. We present a heuristic attack on this protocol that relies on efficient Gröbner bases computations.

The GAP source code referring to the concepts of Chapter 3 and Chapter 4 is provided in Appendix A and Appendix B, respectively.

# Chapter 2

# The Twisted Conjugacy Problem

In this chapter we present a result of Bogopolski, Martino, Maslakova, and Ventura on the solvability of the conjugacy problem in free-by-cyclic groups [6]. To familiarize the reader with the general framework, we introduce free-by-cyclic groups and the twisted conjugacy relation in the first two sections. A brief insight into the subgroup structure of free groups is given in Section 2.3. Equipped with this knowledge, we may proceed with the proof of the main result in Section 2.4. It is shown in [6] that solving the conjugacy problem in free-by-cyclic groups reduces to solving the twisted conjugacy problem in free groups. To this end, two fundamental algorithmic results of Maslakova [22] and Brinkmann [7] are used. Maslakova proved that one can algorithmically find generating sets for the fixed subgroups of free group automorphisms. Brinkmann provided an effective algorithm to determine whether two elements in a free group are mapped to each other by some power of a given automorphism. Finally, we conclude the chapter by giving an overview of some other results known in this context. As already mentioned, this chapter strongly relies on [6].

## 2.1 Free-by-cyclic groups

We embark in this section upon a discussion on free-by-cyclic groups $G$. First we analyze the internal structure of such a group, then we deal with the reverse situation of how to construct $G$ out of certain components. From the preceding we may deduce a presentation for $G$ and a normal form for its elements. The background material in this section is taken from [16] and [18].

**Remark.** Throughout this work, we shall use the notation of [6], so that $w\phi$ denotes the image of an element $w$ in a specified group $G$ under the morphism $\phi$.

**Definition 2.1.** A group $G$ is an *extension* of $N$ by $H$ if $G$ has a normal subgroup $N$ such that the quotient group $G/N$ is isomorphic to $H$. A group $G$ is said to be *free-by-cyclic* if it has a free normal subgroup $F$ such that the quotient group $C \cong G/F$ is cyclic. In other words, $G$ is *free-by-cyclic* if it can be expressed as a *cyclic extension* of $F$ by $C$. If $F$ is finitely generated, then we say that $G$ is a *[f.g.free]-by-cyclic* group.

Let $G$ be a [f.g.free]-by-cyclic group, let $\{x_1, ..., x_n\}$ be a free basis for $F$, and let $C$ be a cyclic group of order $m$. Choose $t$ to be an element of $G$ such that $tF$ is a generator for $G/F$, so that $t^m = h$ for some $h \in F$. The latter follows from $t^m F = (tF)^m = F$ by the group law in the quotient group $G/F$ and since $tF$ is a generator for the cyclic group $G/F$ of order $m$. Suppose that the isomorphism $\psi : C \to G/F$ is given by $c\psi = tF$ where $c$ is a generator for $C$. The map $\theta : C \to Aut(F)$ given by

$$\begin{array}{rcccl} \theta: & C & \to & G/F & \to & Aut(F) \\ & c & \mapsto & tF & \mapsto & x_i\phi := t^{-1}x_i t \quad (1 \leq i \leq n) \end{array}$$

is a homomorphism between $C$ and $Aut(F)$.

We now consider the converse situation of how to recapture $G$ from $F$ and $C$. Clearly, $G$ also involves an action of $C$ on $F$.

Given a finitely generated free group $F$, an element $h \in F$, a positive integer $m$, and the above homomorphism $\theta$ where $\phi \in Aut(F)$ satisfies the conditions

- $\phi^m$ is the inner automorphism of $F$ given by right conjugation by $h$, and
- $\phi$ fixes $h$,

one can reconstruct the [f.g.free]-by-cyclic group $G$ (up to isomorphism). Take $G$ to be the set of ordered pairs $\{(t^k, w) : 0 \leq k \leq m-1, w \in F\}$ equipped with the multiplication rule

$$(t^i, w_1)(t^j, w_2) = \begin{cases} (t^{i+j}, (w_1\phi^j)w_2) & \text{if } i+j < m \\ (t^{i+j-m}, h(w_1\phi^j)w_2) & \text{if } i+j \geq m. \end{cases} \quad (2.1)$$

It can be verified that $G$ is a group with the desired properties. As we will see in Proposition 2.2, every element of $G$ can be uniquely written as $t^k w$ for some $0 \leq k \leq m-1$ and $w \in F$. Here this element is represented by the ordered pair $(t^k, w)$. Multiplying two elements $t^i w_1$ and $t^j w_2$ of $G$, we get

$$(t^i w_1)(t^j w_2) = t^i t^j (t^{-j} w_1 t^j) w_2 = t^{i+j}(w_1 \phi^j) w_2.$$

This leads us in a natural way to the above formula (2.1). If $i + j < m$, then the first case applies. If $i + j \geq m$, then we have $t^{i+j} = t^{i+j-m}h$ which corresponds to the second case.

Taking into account the previous considerations, we deduce the following presentation for $G$

$$M_\phi = \left\langle x_1, ..., x_n, t \mid t^{-1}x_i t = x_i\phi,\ t^m = h \right\rangle,$$

whereas the relation $t^m = h$ is not present when $m = \infty$. In the latter case $M_\phi = F \rtimes_\phi \mathbb{Z}$, that is, $M_\phi$ is the semidirect product of $F$ by $\mathbb{Z} = \langle t \rangle$ with respect to $\phi$.

Up to isomorphism, $M_\phi$ does not depend on $\phi$, but only on the conjugacy class in $Out(F) = Aut(F)/Inn(F)$ of its induced outer automorphism. We refer the reader to Lemma 2.1 of [4] for a proof of this statement in the case where $m = \infty$.

It is conjectured that most one-relator groups are [f.g.free]-by-cyclic. However, this has not yet been expressed in a precise form at this time.

The next proposition states that we have a normal form for elements in $M_\phi$, which is algorithmically computable from a given arbitrary word on the generators.

**Proposition 2.2.** *Every element in $M_\phi$ can be uniquely represented as a word of the form $t^r u$, where $r$ is an integer and $u \in F$. Moreover, when $C$ is finite of order $m$, we may assume that $0 \leq r \leq m - 1$.*

*Proof.* Using the defining relations of $M_\phi$ under the form $wt = t(w\phi)$ and $wt^{-1} = t^{-1}(w\phi^{-1})$ (the inverse of $\phi$ is given by right conjugation by $t^{-1}$) for $w \in F$, we can push all the letters $t$ to the left, in every element of $M_\phi$. This establishes the first claim.

When $C$ is finite of order $m$ (that is, when the relation $t^m = h\ (\in F)$ is present), we may write $t^n w = t^r t^{qm} w = t^r h^q w = t^r w'$ for $n = qm + r$ with $q \in \mathbb{Z}$ and $0 \leq r \leq m - 1$ (where $q$ and $r$ are given by the Euclidean algorithm). □

## 2.2 Twisted conjugacy

We start with the definition of the twisted conjugacy relation.

**Definition 2.3.** Let $G$ be an arbitrary group and let $\phi \in Aut(G)$. Two elements $u, v \in G$ are said to be *$\phi$-twisted conjugated*, denoted by $u \sim_\phi v$, if there exists an element $g \in G$ with $(g\phi)^{-1} u g = v$.

**Proposition 2.4.** *The relation $\sim_\phi$ is an equivalence relation.*

*Proof.* For all $u, v, w \in G$, we verify the three requirements on $\sim_\phi$.

- It follows from Definition 2.3 that $u \sim_\phi u$ by taking $g = 1$. Hence, $\sim_\phi$ is reflexive.

- If $u \sim_\phi v$, then there exists $g \in G$ such that $(g\phi)^{-1}ug = v$. This implies that $(g\phi)vg^{-1} = u$, so that $v \sim_\phi u$ with $\phi$-twisted conjugating element $g^{-1}$. Thus, $\sim_\phi$ is symmetric.

- If $u \sim_\phi v$ and $v \sim_\phi w$, then there exist $g, h \in G$ such that $(g\phi)^{-1}ug = v$ and $(h\phi)^{-1}vh = w$. So

$$\begin{aligned} w &= (h\phi)^{-1}(g\phi)^{-1}ugh \\ &= ((gh)\phi)^{-1}ugh \end{aligned}$$

showing that $u \sim_\phi w$ with $\phi$-twisted conjugating element $gh$. This checks that $\sim_\phi$ is transitive.

$\square$

The equivalence class of $u \in G$ under $\sim_\phi$ is known as the orbit of $u$ and will be given by $\mathrm{TwOrb}(u) = \{v \in G \mid v \sim_\phi u\}$. It is also called the $\phi$-twisted conjugacy class of $u$ (or Reidemeister class). A general algorithm for recognizing $\phi$-twisted conjugacy classes is not yet known. However, in Theorem 2.16 such an algorithm is provided in the case where $\phi$ is an automorphism of a finitely generated free group.

We next define the twisted conjugacy problem. This decision problem will play a central role in Section 2.4.

- The $\phi$-*twisted conjugacy problem* for a given $\phi \in Aut(G)$ is the algorithmic problem of deciding, for two arbitrary words $u, v \in G$, whether or not $u \sim_\phi v$.

- The *twisted conjugacy problem* for a group $G$ is the decision problem of determining whether the $\phi$-twisted conjugacy problem is solvable for any $\phi \in Aut(G)$.

**Remark.** Twisted conjugacy can be considered as a generalization of the ordinary conjugacy relation in groups. In particular, the *id*-twisted conjugacy problem for a group $G$ corresponds to the standard conjugacy problem in $G$.

Let $\phi$ be an automorphism of a free group $F$ and let $u \in F$. Then we denote by $\mathrm{Orb}(u) = \{v \in F \mid v = \phi^m(u) \text{ for some } m \in \mathbb{Z}\}$ the $\phi$-orbit of $u$. The following lemma will be used in the proof of Proposition 2.15.

**Lemma 2.5.** *Any $\phi$-twisted conjugacy class in $F$ is a union of $\phi$-orbits.*

*Proof.* We have to show that if $v, w \in \mathrm{Orb}(u)$, for two words $v, w \in F$, then $v, w \in \mathrm{TwOrb}(u)$. The converse of this statement is not necessarily true (see example below). First, notice that $u \sim_\phi u\phi$ since $u = (u\phi)^{-1}(u\phi)u$. Also, $u\phi^k \sim_\phi u\phi^{k+1}$ $\forall u \in F$ and $\forall k \in \mathbb{Z}$. This follows from the fact that $u\phi^{k+1} = ((u^{-1})\phi^{k+1})^{-1}(u\phi^k)(u^{-1}\phi^k)$. By transitivity of the equivalence relation $\sim_\phi$, we get that $\mathrm{Orb}(u) \subset \mathrm{TwOrb}(u)$.

Now, take $v \in F$ such that $v \in \mathrm{TwOrb}(u)$, but $v \notin \mathrm{Orb}(u)$. We again obtain that $\mathrm{Orb}(v) \subset \mathrm{TwOrb}(u)$, by transitivity of $\sim_\phi$. Repeating this process yields the desired result, that is,

$$\mathrm{TwOrb}(u) = \cup_{i=0}^\infty \mathrm{Orb}(u_i),$$

where $u_i \notin \mathrm{Orb}(u_j)$ $\forall\, i \neq j$ and $u_i \in \mathrm{TwOrb}(u)$ $\forall i \in \mathbb{N}$. $\square$

**Example.** Let us consider the automorphism of the free group $F_2$ given by the rule
$$\phi : \begin{cases} y_1 & \longmapsto & y_2 \\ y_2 & \longmapsto & y_1 y_2, \end{cases}$$
where $y_1$ and $y_2$ are the free generators of $F_2$.
Furthermore, let $u := y_2 y_1$ and $v := y_1^2$ be words in $F_2$. Then $u \sim_\phi v$ since

$$(y_1\phi)^{-1} y_2 y_1 y_1 = y_2^{-1} y_2 y_1^2 = y_1^2.$$

Hence, $v \in \mathrm{TwOrb}(u)$. However, we show next that $v \notin \mathrm{Orb}(u)$.
The inverse of the automorphism $\phi$ is given by the rule

$$\phi^{-1} : \begin{cases} y_1 & \mapsto & y_2 y_1^{-1} \\ y_2 & \mapsto & y_1 \end{cases}$$

The images of $v$ under positive powers of $\phi$ can be computed without performing any cancellations:

$$y_2 y_1 \mapsto y_1 y_2 y_2 \mapsto y_2 y_1 y_2 y_1 y_2 \mapsto y_1 y_2 y_2 y_1 y_2 y_2 y_1 y_2 \mapsto \ldots.$$

Similarly, under successive applications of $\phi^{-1}$ we get

$$y_2 y_1 \mapsto y_1 y_2 y_1^{-1} \mapsto y_2 y_1^{-1} \cdot y_1 y_1 y_2^{-1} \mapsto y_1 y_2 y_1^{-1} y_1^{-1} \mapsto y_2 y_1 y_2^{-1} y_1 y_2^{-1} \mapsto \ldots.$$

Here the first cancellation appears in $\phi^{-2}(y_2 y_1)$ while applying $\phi^{-1}$ on $y_1 y_2$ in the previous step. Because of these cancellations, the length of $y_2 y_1$ under negative powers of $\phi$ cannot be controlled as easily as before. But the form of $\phi^{-k}$ for $k \geq 0$ can be foreseen and will never be $y_1^2$. Therefore, $v \notin \mathrm{Orb}(u)$. This example suggests that there exist words $u$ and $v$ in $F_2$ such that $v \in \mathrm{TwOrb}(u)$, but $v \notin \mathrm{Orb}(u)$ (as stated in the proof of Lemma 2.5).

## 2.3 Graph representation for subgroups of free groups

To solve the twisted conjugacy problem in free groups, we use an algorithm for deciding membership in a subgroup of a free group. To this end, we first associate a graph to such a subgroup. This requires the topological notion of a folding of graphs, introduced by Stallings in [32]. A more combinatorial view on this topic is given in [20], on which the following is based.

First, we need some basic definitions. Let $X = \{x_1, ..., x_n\}$ be a finite alphabet and let $\Sigma = X \cup X^{-1}$.

**Definition 2.6.** An *X-digraph* $\Gamma$ is a tuple consisting of a nonempty set of vertices $X^0$, a set of edges $X^1$ and three mappings $\alpha : X^1 \to X^0$, $\omega : X^1 \to X^0$ and $^- : X^1 \longrightarrow X^1$ (meaning the beginning, the end and the inverse of an edge) such that $\bar{\bar{e}} = e$, $\bar{e} \neq e$ and $\alpha(e) = \omega(\bar{e})$ $\forall e \in X^1$. Furthermore, every edge $e$ is labeled by a letter from $X$, denoted by $\mu(e)$, such that $\mu(\bar{e}) = (\mu(e))^{-1}$.

**Definition 2.7.** Let $\Gamma$ be an $X$-digraph. We say that $\Gamma$ is *folded* if for each vertex $v$ of $\Gamma$ and each letter $a \in X$ there is at most one edge in $\Gamma$ with origin $v$ and label $a$ and there is at most one edge with terminus $v$ and label $a$.

Suppose $\Gamma$ is an $X$-digraph. If it happens that two edges have the same initial vertex and the same label, then we may identify these edges and their terminal vertices, and give the resulting edge the same label. This operation is known as a *folding*.

**Definition 2.8.** Let $\Gamma$ be an $X$-digraph and let $v$ be a vertex of $\Gamma$. *The language of $\Gamma$ with respect to $v$* is defined as

$$L(\Gamma, v) = \{\mu(p) \mid p \text{ is a reduced path in } \Gamma \text{ from } v \text{ to } v\}$$

If a word $w \in L(\Gamma, v)$, then we say that $w$ is *accepted* by $(\Gamma, v)$.

**Remark.** The words in $L(\Gamma, v)$ are not necessarily freely reduced. However, if $\Gamma$ is a folded $X$-digraph, then all the words in $L(\Gamma, v)$ are freely reduced.

**Definition 2.9.** The *core of an X-digraph* $\Gamma$ *at a vertex* $v$ of $\Gamma$ is the (connected) subgraph of $\Gamma$ containing all reduced paths from $v$ to $v$.

The proof of the following proposition is omitted, since it is quite cumbersome (see, for example, [20]).

**Proposition 2.10.** *For each subgroup $H$ of a free group $F$ on $X$, there exists a unique $X$-digraph $(\Gamma(H), 1_H)$, up to isomorphism of $X$-digraphs, such that the graph $\Gamma(H)$ is folded, connected and is a core graph with respect to $1_H$. Moreover, the language of $\Gamma(H)$ with respect to $1_H$ is $H$, that is $L(\Gamma(H), 1_H) = H$.*

**Remark.** If $H$ is finitely generated, then this graph is finite.

Let $F$ be a free group on $X$ and let $H = \langle h_1, ..., h_n \rangle$ a finitely generated subgroup of $F$ for some words $h_1, ..., h_n$ in $\Sigma = X \cup X^{-1}$.

**Proposition 2.11.** *There is an algorithm constructing the graph $\Gamma(H)$.*

*Proof.* The idea is to start with a wedge of $n$ circles with the words $h_1, ..., h_n$ written on them and then perform all possible foldings. The resulting $X$-digraph will be $\Gamma(H)$. Notice that, starting with a graph having $k := |h_1| + ... + |h_n|$ edges, we obtain the graph $\Gamma(H)$ in at most $k$ steps (since each folding reduces the number of edges by one). We refer the reader to [3], [20] or [32] for further details. □

To point out the usefulness of the graph $\Gamma(H)$, we present here an algorithm for detecting membership in $H$. It is worth mentioning that algorithms involving $\Gamma(H)$ are quite efficient since they usually run in polynomial time.

**Proposition 2.12.** *The membership problem for $H$ in $F$ is solvable. That is, there exists an algorithm deciding for an arbitrary word $w$ in $F$ whether or not $w$ lies in the subgroup $H$.*

*Proof.* First, we construct the graph $\Gamma = \Gamma(H)$ (as done in Proposition 2.11). Recall that $\Gamma$ is folded and that $L(\Gamma, 1_H) = H$.

Given a word $w = x_1 \cdots x_m \in F$ with $x_i \in \Sigma$, we check whether $\Gamma(H)$ accepts $w$ as follows. If there is an edge with origin $1_H$ and label $x_1$ (such an edge is unique if it exists, since $\Gamma$ is folded), then we move to the terminal vertex of this edge, denoted by $v_1$. If not, we terminate the process. Then, we verify whether there is an edge with origin $v_1$ and label $x_2$. Repeating this procedure at most $m$ times, we conclude that $w \notin L(\Gamma, 1_H)$ if $v_m \neq 1_H$ or if the process has already terminated for some $k < m$ and we conclude that $w \in L(\Gamma, 1_H)$ if $v_m = 1_H$. □

## 2.4 The conjugacy problem for free-by-cyclic groups

We only state the results of Maslakova [22] and Brinkmann [7] that will be used later on; their proof is being omitted, since they are quite involved. They both make use of the theory of train tracks developed by Bestvina and Handel in [2].

**Theorem 2.13** (Maslakova, 2003)**.** *There exists an algorithm to compute a finite generating set for the fixed point subgroup of an arbitrary automorphism of a free group of finite rank.*

**Theorem 2.14** (Brinkmann, 2008). *Let $\phi$ be an automorphism of a finitely generated free group $F$. There exists an explicit algorithm that, given two elements $u, v \in F$, decides whether there exists an integer $k$ such that $u\phi^k$ is conjugate to $v$. If such an integer $k$ exists, then the algorithm will compute $k$ as well.*

Recall that we have a normal form in $M_\phi$ whereby every element can be written in a unique way as $t^r u$, where $r$ is an integer and $u \in F$. When $C$ is of order $m < \infty$, we may assume that $0 \leq r \leq m - 1$.
Conjugating $t^r u$ by an arbitrary element $t^k g$, we get

$$\begin{aligned}(t^k g)^{-1}(t^r u)(t^k g) &= g^{-1}t^{-k}t^r u t^k g = g^{-1}t^r t^{-k} u t^k g \\ &= g^{-1}t^r(u\phi^k)g = t^r t^{-r} g^{-1} t^r (u\phi^k) g \\ &= t^r (g\phi^r)^{-1}(u\phi^k)g.\end{aligned}$$

Therefore, two elements $t^r u$ and $t^s v$ in $M_\phi$ (with $0 \leq r, s \leq m - 1$ in the case where $|C| = m < \infty$ and $u, v$ in F) are conjugate in $M_\phi$ if and only if $r = s$ and $v = (g\phi^r)^{-1}(u\phi^k)g$ for some integer $k$, that is, $v \sim_{\phi^r} u\phi^k$. This is the key fact in the proof of the following proposition.

**Proposition 2.15.** *Let $F$ be a finitely generated free group. If the twisted conjugacy problem is solvable in $F$, then the standard conjugacy problem is solvable in $M_\phi$, for every $\phi \in Aut(F)$.*

*Proof.* Given two elements $t^r u$ and $t^s v$ in $M_\phi$, we have to decide whether they are conjugate in $M_\phi$ and find a conjugating element if it exists.

*Case $r = 0$.* Notice that $u$ can only be conjugate in $M_\phi$ to other elements $v$ of the base group $F$ (since $r = s = 0$). Moreover, $u$ is conjugate to $v$ in $M_\phi$ if and only if $u\phi^k$ is conjugate to $v$ for some integer $k$. Since this is decidable by Theorem 2.14, we can decide whether $u, v \in F$ are conjugate in $M_\phi$.

*Case $r \neq 0$.* By Lemma 2.5, $u\phi^k \sim_{\phi^r} u\phi^{k \pm r}$. Hence, $t^r u$ and $t^s v$ are conjugate in $M_\phi$ if and only if $r = s$ and $v = (g\phi^r)^{-1}(u\phi^k)g$ for some integer $0 \leq k \leq |r| - 1$. Thus, finitely many checks of the twisted conjugacy problem in $F$ provide a solution to the standard conjugacy problem in $M_\phi$. □

**Theorem 2.16.** *Let $F$ be a finitely generated free group. The twisted conjugacy problem is solvable in $F$.*

*Proof.* Let $\phi \in Aut(F)$ and suppose that $u, v \in F$ are given. We need to algorithmically decide whether $u \sim_\phi v$.

Choose a free basis for $F$. Adding a new letter $z$, we get a free basis for $F' = F * \langle z \rangle$. Let $\phi' \in Aut(F')$ be the extension of $\phi$ defined by $z\phi' = uzu^{-1}$ and let $\gamma_y$ denote the inner automorphism of $F'$ given by right conjugation by $y \in F'$, $x\gamma_y = y^{-1}xy$.

*Claim.* It holds that $u \sim_\phi v$ if and only if $Fix(\phi'\gamma_v)$ contains an element of the form $g^{-1}zg$ for some $g \in F$ ($g$ is then a valid $\phi$-twisted conjugating element).

Suppose that there exists an element $g \in F$ such that $(g\phi)^{-1}ug = v$ (that is $u \sim_\phi v$). Then,

$$\begin{aligned}(g^{-1}zg)\phi'\gamma_v &= v^{-1}(g\phi)^{-1}uzu^{-1}(g\phi)v \\ &= g^{-1}u^{-1}(g\phi)(g\phi)^{-1}uzu^{-1}(g\phi)(g\phi)^{-1}ug \\ &= g^{-1}zg.\end{aligned}$$

Conversely, if $g^{-1}zg \in Fix(\phi'\gamma_v)$ for some $g \in F$, then

$$g^{-1}zg = (g^{-1}zg)\phi'\gamma_v = v^{-1}(g\phi)^{-1}uzu^{-1}(g\phi)v.$$

Hence, $gv^{-1}(g\phi)^{-1}uz = zgv^{-1}(g\phi)^{-1}u$ which means that $gv^{-1}(g\phi)^{-1}u$ (containing no occurrences of $z$) commutes with $z$. This implies that $gv^{-1}(g\phi)^{-1}u = 1$. In other words, $(g\phi)^{-1}ug = v$ so that $u \sim_\phi v$ with $g$ being a $\phi$-twisted conjugating element.

By Theorem 2.13, we can compute a finite generating set for the subgroup $Fix(\phi'\gamma_v)$. Therefore, it is algorithmically decidable whether $Fix(\phi'\gamma_v)$ contains an element of the form $g^{-1}zg$ for some $g \in F$ (by modifying the algorithm in the proof of Proposition 2.12). One can, for example, check if the finite core-graph of $Fix(\phi'\gamma_v)$ contains a loop labeled $z$ at some vertex connected to the base-point by a path whose label does not use the letter $z$. If this is the case, the label of such a path provides the $\phi$-twisted conjugating element. □

Finally, combining Proposition 2.15 with Theorem 2.16, we immediately obtain the following result.

**Theorem 2.17** (Bogopolski-Martino-Maslakova-Ventura, 2006). *The conjugacy problem in [f.g.free]-by-cyclic groups is solvable.*

## 2.5 Some related results

The result of Theorem 2.17 was already known for some particular cases, which we describe in Subsection 2.5.1. Then, in Subsections 2.5.2 and 2.5.3, we explain how this result was adjusted to free-by-free and [free abelian]-by-free groups by Bogopolski, Martino and Ventura [5]. The notion of orbit decidability and how it relates to the solvability of the conjugacy problem in free-by-free groups are the subject of Subsection 2.5.2. In Subsection 2.5.3 we shall be mainly concerned with the construction of an explicit orbit undecidable subgroup of $Aut(\mathbb{Z}^n) = GL_n(\mathbb{Z})$ ($n \geq 4$) giving rise to a [free abelian]-by-free group with unsolvable conjugacy problem. The material for the last two subsections is taken from [5].

### 2.5.1 From free-by-cyclic to hyperbolic groups

In this subsection we use the notation as introduced in Section 2.1.

An automorphism $\phi$ of the finitely generated free group $F$ is said to have no nontrivial periodic conjugacy classes if one cannot find an integer $k \neq 0$ and elements $g, h \in F$ such that $1 \neq g\phi^k = h^{-1}gh$. If the cyclic group $C$ is infinite, this is equivalent to the fact that the [f.g.free]-by-cyclic group $M_\phi$ has no $\mathbb{Z} \oplus \mathbb{Z}$ subgroups. Therefore, $M_\phi$ is hyperbolic (in the sense of Gromov) which has been proved by Brinkmann in [8]. Furthermore, when $C$ is finite, the group $M_\phi$ is virtually free (that is, having a finite-index subgroup which is free), thus $M_\phi$ is quasi-isometric to a free group, and hence hyperbolic (see [15]). Since for hyperbolic groups the conjugacy problem is solvable (see, for example, Lemma 8 of [15]), the same holds for $M_\phi$ in the previous cases. To summarize, for groups which are both [f.g.free]-by-cyclic and hyperbolic, it is well known that they have solvable conjugacy problem. So this is a special case of Theorem 2.17. But, not all [f.g.free]-by-cyclic groups are hyperbolic and vice-versa.

Suppose that $\phi$ is a virtually inner automorphism, that is $\phi^m \in Inn(F)$ for some non-zero integer $m$. Bardakov, Bokut, and Vesnin proved in [1] that there exists an algorithm to decide, for two words $tu, tv \in M_\phi$ (with $u, v \in F$), whether or not they are conjugated by some element in $F$; this is equivalent to the solvability of the $\phi$-twisted conjugacy problem in $F$. They also showed that, if $\phi \in Aut(F)$ has no nontrivial periodic conjugacy classes, then the $\phi$-twisted conjugacy problem in $F$ is solvable. This mainly relies on the fact that the conjugacy problem for hyperbolic groups is solvable (refer to the previous paragraph).

### 2.5.2 The conjugacy problem for free-by-free groups

Theorem 2.17 is no longer true for [f.g.free]-by-[f.g.free] groups. An explicit construction of such a group for which the conjugacy problem is unsolvable is given by Miller in [25]. An alternative proof of the unsolvability of the conjugacy problem for Miller's group via Mihailova's result is given in [5] .

The *orbit decidability problem* for a subgroup $A \leq Aut(G)$ is the algorithmic problem of deciding for two arbitrary elements $u, v \in G$, whether or not there exists a $\phi \in A$ such that $u\phi$ and $v$ are conjugate to each other in $G$. In [5], Bogopolski, Martino and Ventura showed that orbit decidability is the only obstruction to the solvability of the conjugacy problem in [f.g.free]-by-[f.g.free] groups (in fact they proved this for an even larger family of groups). So the proof of Theorem 2.17 can be extended to this bigger family of groups, except for the step where Brinkmann's result comes into play. Here the stronger problem of orbit decidability arises. This is formally expressed in the next theorem.

**Theorem 2.18.** *Let $F_n$ be the free group on $\{x_1, ..., x_n\}$ ($n \geq 2$), and let $\phi_1, ..., \phi_m$ be automorphisms of $F_n$. Then the [f.g.free]-by-[f.g.free] group*

$$F_n \rtimes_{\phi_1, ..., \phi_m} F_m = \left\langle x_1, ..., x_n, t_1, ..., t_m \mid t_j^{-1} x_i t_j = x_i \phi_j \right\rangle$$

*has solvable conjugacy problem if and only if the subgroup $\langle \phi_1, ..., \phi_m \rangle \leq Aut(F)$ is orbit decidable.*

In the [f.g.free]-by-cyclic case, orbit decidability is equivalent to cyclic subgroups of $Aut(F)$ being orbit decidable, which has already been solved by Brinkmann in [7]. Theorem 2.17 follows as an immediate consequence.

Whitehead's theorem (see [37]) states that the full automorphism group $Aut(F_n)$ of a finitely generated free group $F_n$ is orbit decidable. In light of Theorem 2.18, this result implies that if $\langle \phi_1, ..., \phi_m \rangle = Aut(F)$, then the group $F_n \rtimes_{\phi_1, ..., \phi_m} F_m$ has solvable conjugacy problem.

### 2.5.3 The conjugacy problem for [free abelian]-by-free groups

The analog of Theorem 2.18 for [f.g.free abelian]-by-[f.g.free] groups is as follows.

**Theorem 2.19.** *Let $\mathbb{Z}^n = \langle x_1, ..., x_n \mid [x_i, x_j] \rangle$ be the free abelian group of rank $n \geq 2$, and let $M_1, ..., M_m \in Aut(\mathbb{Z}^n) = GL_n(\mathbb{Z})$. Then the [f.g.free abelian]-by-[f.g.free] group*

$$\mathbb{Z}^n \rtimes_{M_1, ..., M_m} F_m = \left\langle x_1, ..., x_n, t_1, ..., t_m \mid t_j^{-1} x_i t_j = x_i M_j, \ [x_i, x_j] \right\rangle$$

*has solvable conjugacy problem if and only if the subgroup $\langle M_1, ..., M_m \rangle \leq GL_n(\mathbb{Z})$ is orbit decidable.*

Using facts from linear algebra, similar results as in the previous subsection can be obtained to detect [free abelian]-by-free groups with solvable conjugacy problem.

Before giving an explicit example of a [free abelian]-by-free group with unsolvable conjugacy problem, we review a general source of orbit undecidability that will apply not exclusively to free and free abelian groups $G$.

First we need to define a few concepts. Let $G$ be a group. The *stabilizer* of a subgroup $K \leq G$ is given by

$$Stab(K) = \{\phi \in Aut(G) \mid k\phi = k \ \forall k \in K\} \leq Aut(G).$$

The *conjugacy stabilizer* of the subgroup $K$, denoted by $Stab^*(K)$, is the set of automorphisms acting as conjugation on $K$, that is, $Stab^*(K) = Stab(K)Inn(G) \leq Aut(G)$.

**Proposition 2.20.** *Let $G$ be a group, and let $A \leq B \leq Aut(G)$ and $v \in G$ be such that $B \cap Stab^*(\langle v \rangle) = \{Id\}$. If $A \leq Aut(G)$ is orbit decidable, then the membership problem for $A$ in $B$ (that is, given an element $\phi \in B$ to decide whether or not $\phi \in A$) is solvable.*

*Proof.* Given $\psi \in B \leq Aut(G)$, we have to decide whether or not $\psi \in A$. Let $w = v\psi \in G$. Observe that

$$\{\phi \in B \mid v\phi \sim w\} = B \cap (Stab^*(\langle v \rangle)\psi) = (B \cap Stab^*(\langle v \rangle))\psi = \{\psi\}.$$

Then there exists $\varphi \in A \cap \{\phi \in B \mid v\phi \sim w\}$ if and only if $\psi \in A$. Orbit decidability of the subgroup $A \leq Aut(G)$ allows to solve the membership problem for $A$ in $B$. □

It is well known that for a free group $F$ the free product $F * F$ is again a free group and so the membership problem for finitely generated subgroups in $F * F$ is solvable. However, in [24], Mihailova showed that solvability of the membership problem is not preserved by direct products of free groups $F \times F$. For a detailed treatment of the proof, we refer the reader to [25, 26].

**Theorem 2.21** (Mihailova, 1959). *Let $F_n$ be a finitely generated free group of rank $n \geq 2$. Then the group $G = F_n \times F_n$ has a finitely generated subgroup $L$ such that the membership problem for $L$ in $G$ is recursively unsolvable.*

The proof of Mihailova's theorem is based on the following construction. Start with a finitely presented group $H = \langle s_1, ..., s_n \mid R_1, ..., R_m \rangle$ and consider the subgroup $A = \{(x, y) \in F_n \times F_n \mid x =_H y\} \leq F_n \times F_n$. It is not difficult to verify that $A = \langle (1, R_1), ..., (1, R_m), (s_1, s_1), ..., (s_n, s_n) \rangle$. Clearly, if the word problem for $H$ is unsolvable, the membership problem for $A$ in $F_n \times F_n$ is unsolvable too.

We may focus our attention on two generator groups since by the Higman-Neumann-Neumann embedding theorem (see, for example, [26]) every countable group can be embedded in a two generator group with the same number of relations. Together with Mihailova's result and Proposition 2.20, we obtain the following theorem.

**Theorem 2.22.** *Let $G$ be a finitely generated group such that $F_2 \times F_2$ embeds in $Aut(G)$ in such a way that the image intersects trivially with $Stab^*(\langle v \rangle)$ for some $v \in G$. Then $Aut(G)$ contains an orbit undecidable subgroup. In other words, there exist $G$-by-[f.g.free] groups with unsolvable conjugacy problem.*

Let us concentrate on free abelian groups $G = \mathbb{Z}^n$ ($n \geq 4$) and construct orbit undecidable subgroups of $GL_n(\mathbb{Z})$ that will correspond to the first known examples of [f.g.free abelian]-by-[f.g.free] groups with unsolvable conjugacy problem.

**Proposition 2.23.** *For $n \geq 4$, $GL_n(\mathbb{Z})$ contains finitely generated orbit undecidable subgroups.*

*Proof.* We only outline the key ideas of the proof (for further details see [5]). First, notice that $F_2$ embeds in $GL_2(\mathbb{Z})$ so that $F_2 \times F_2$ embeds in $GL_4(\mathbb{Z})$.

*Step 1.* Consider the free subgroup $\langle P, Q \rangle \leq GL_2(\mathbb{Z})$, where $P = \begin{pmatrix} 1 & 1 \\ 1 & 2 \end{pmatrix}$ and $Q = \begin{pmatrix} 2 & 1 \\ 1 & 1 \end{pmatrix}$. (Note that $F_2 \cong \langle P, Q \rangle \leq SL_2(\mathbb{Z})$.)

One has to check that

1. $Stab^*(\langle (1,0) \rangle) = Stab(\langle (1,0) \rangle) = \left\{ \begin{pmatrix} 1 & 0 \\ n & \pm 1 \end{pmatrix} \mid n \in \mathbb{Z} \right\}$, and

2. $\langle P, Q \rangle \cap Stab(\langle (1,0) \rangle) = \left\langle \begin{pmatrix} 1 & 0 \\ 12 & 1 \end{pmatrix} \right\rangle$.

*Step 2.* Choose a free subgroup of rank 2 $\langle P', Q' \rangle \leq \langle P, Q \rangle$ such that $\langle P', Q' \rangle \cap Stab(\langle (1,0) \rangle) = \{Id\}$. Consider, for $n \geq 4$,

$$B = \left\langle \left( \begin{array}{c|c} P' & 0 \\ \hline 0 & Id \end{array} \right), \left( \begin{array}{c|c} Q' & 0 \\ \hline 0 & Id \end{array} \right) \left( \begin{array}{c|c} Id & 0 \\ \hline 0 & P' \end{array} \right) \left( \begin{array}{c|c} Id & 0 \\ \hline 0 & Q' \end{array} \right) \right\rangle \leq GL_n(\mathbb{Z}).$$

Clearly, $B \cong F_2 \times F_2$. By construction, $B$ intersects trivially with the stabilizer of $v = (1, 0, 1, 0, ..., 0) \in \mathbb{Z}^n$.

*Step 3.* Find a finitely generated subgroup $A \leq B$ with unsolvable membership problem in $B$ via Mihailova's construction. Applying Proposition 2.20 to $G = \mathbb{Z}^n$, the subgroup $A \leq GL_n(\mathbb{Z}) = Aut(\mathbb{Z}^n)$ reveals to be orbit undecidable. □

**Corollary 2.24.** *There exist $\mathbb{Z}^n$-by-[f.g.free] groups with unsolvable conjugacy problem $(n \geq 4)$.*

# Chapter 3

# The Doubly-Twisted Conjugacy Problem

In this chapter we shall discuss recent work by Staecker on doubly-twisted conjugacy relations (see [29] and [31]). Section 3.1 establishes a link between doubly-twisted conjugacy, the equalizer subgroup and Post's Correspondence Problem. Then, in Section 3.2, the remnant property for homomorphisms is introduced. As we shall see in Sections 3.3 and 3.4, this property is extremely useful to provide criteria for membership of distinct words in different doubly-twisted conjugacy classes as well as for the equalizer subgroup to be trivial. The last section is devoted to some open problems related to the equalizer subgroup.

## 3.1 Doubly-twisted conjugacy and the equalizer subgroup

Let $G, H$ be finitely generated free groups and let $\phi, \psi : G \to H$ be homomorphisms.

**Definition 3.1.** Two elements $u, v \in H$ are said to be *doubly-twisted conjugated* if there exists an element $g \in G$ such that $v = (g\phi)^{-1} u (g\psi)$.

Doubly-twisted conjugacy is an equivalence relation and we denote by $[u]$ the equivalence class of $u \in H$. No general algorithm exists to decide doubly-twisted conjugacy in free groups.

The *doubly-twisted conjugacy problem* is the following decision problem. Given a pair of homomorphisms $\phi, \psi$ between two arbitrary groups $G$ and $H$ and a pair of elements $u, v \in G$, decide whether or not there is an element $g \in G$ such that $v = (g\phi)^{-1} u (g\psi)$.

**Remark.** Doubly-twisted conjugacy is a generalization of the twisted conjugacy relation defined in the previous chapter (simply set $\psi = id$).

**Definition 3.2.** Given two homomorphisms $\phi, \psi : G \to H$, the *equalizer subgroup* $Eq(\phi, \psi) \leq G$ is defined as

$$Eq(\phi, \psi) = \{g \in G \mid g\phi = g\psi\}. \tag{3.1}$$

Let $\langle z \rangle$ be the free group generated by $z$. Let $\widehat{G} = G * \langle z \rangle$ and $\widehat{H} = H * \langle z \rangle$, where $*$ denotes the free product. Furthermore, let $\phi^v$ be given by $\phi$ followed by right conjugation by $v \in H$, that is, $g\phi^v = v^{-1}(g\phi)v$ for $g \in G$. The following Lemma holds for any groups $G$ and $H$ (not necessarily free).

**Lemma 3.3.** *Let $\widehat{\phi} : \widehat{G} \to \widehat{H}$ be the extension of $\phi$ given by $z\widehat{\phi} = uzu^{-1}$ for a given $u \in H$ and let $\widehat{\psi} : \widehat{G} \to \widehat{H}$ be the extension of $\psi$ given by $z\widehat{\psi} = z$. Then $[u] = [v]$ if and only if there is some $g \in G$ with $g^{-1}zg \in Eq(\widehat{\phi}^v, \widehat{\psi})$.*

*Proof.* The proof is very similar to what we have already seen in the proof of Theorem 2.16. □

Goldstein and Turner [12] proved that the subgroup $Eq(\phi, id) = Fix(\phi) = \{w \in F \mid w\phi = w\}$ is finitely generated for any endomorphism $\phi$ of a finitely generated free group $F$. Moreover, Imrich and Turner [17] showed that the rank of $Fix(\phi)$ cannot exceed the rank of $F$ (by generalizing the same result for automorphisms due to Bestvina and Handel [2]).

So far very few algorithmic results on fixed subgroups of free groups are known. In 2003, Maslakova [22] introduced a complicated algorithm to compute a basis for $Fix(\phi)$ if $\phi$ is an automorphism of $F$. In this case, one can easily decide whether $Fix(\phi)$ contains a certain element of $F$ using Stallings graph, and thus decide whether two words belong to the same twisted conjugacy class (see proof of Theorem 2.16). However, the corresponding problems for endomorphisms of $F$ are still open at the time of writing.

The equalizer subgroup of two homomorphisms is not always finitely generated (only in case one of them is injective which was shown by Goldstein and Turner in [12]). A general approach to compute a basis for such a subgroup seems to be out of reach. Also, an algorithm to verify whether or not a given word $w$ belongs to $Eq(\phi, \psi)$ is not yet known. Since computing the equalizer subgroup of homomorphisms reveals to be difficult, Lemma 3.3 cannot be applied to recognize doubly-twisted conjugacy classes.

### 3.1.1 Relation to PCP

Post's Correspondence Problem (PCP) is one of the most impressive undecidable problems in computer science. We formulate the same problem for free groups. Here the setting is slightly different because of the cancellations that might occur between consecutive words.

PCP in free groups can be formulated as follows. Given two homomorphisms $\phi$ and $\psi$ between finitely generated free groups $G$ and $H$ of ranks $n$ and $m > 1$ respectively, is there an element $1 \neq w \in G$ such that $w\phi = w\psi$?

More precisely, if $\{a_1, ..., a_n\}$ is a generating set for $G$, then consider the two $n$-tuples of elements of $H$, namely $\{a_1\phi, ..., a_n\phi\}$ and $\{a_1\psi, ..., a_n\psi\}$. Is there a correspondence between them? In other words, is there a finite sequence of indices $i_1, ..., i_k \in \{1, ..., n\}$, $k \geq 1$, such that

$$(a_{i_1}\phi)\cdots(a_{i_k}\phi) = (a_{i_1}\psi)\cdots(a_{i_k}\psi)?$$

This problem is equivalent to the problem of deciding whether the equalizer subgroup $Eq(\phi, \psi)$ is trivial or not.

We now focus our attention on the following question. Given $Eq(\phi, \psi)$, what information can we obtain about $Eq(\widehat{\phi^v}, \widehat{\psi})$, and vice versa? By Lemma 3.3, this could in turn yield a relation between membership in the equalizer subgroup $Eq(\phi, \psi)$ and deciding doubly-twisted conjugacy relations.

**Proposition 3.4.** *Let $w, v \in G$ such that $w\phi \in \langle v \rangle$. Then $w \in Eq(\phi, \psi)$ if and only if $w \in Eq(\widehat{\phi^v}, \widehat{\psi})$.*

*Proof.* Let $w\phi = v^m$ for some $m \in \mathbb{Z}$. First assume that $w \in Eq(\phi, \psi)$. Then

$$w\widehat{\phi^v} = v^{-1}(w\phi)v = v^{-1}v^m v = v^m = w\phi = w\psi = w\widehat{\psi}.$$

Conversely, if $w \in Eq(\widehat{\phi^v}, \widehat{\psi})$, then

$$w\psi = w\widehat{\psi} = w\widehat{\phi^v} = v^{-1}(w\phi)v = v^{-1}v^m v = v^m = w\phi.$$

□

## 3.2 Homomorphisms with remnant

Informally, a homomorphism $\phi$ has remnant if the images of the generators under $\phi$ have limited cancellation when multiplied together. We give here a precise definition and an example of this concept which was first introduced by Wagner in [36].

**Definition 3.5.** Let $G$ be a finitely generated free group on the generators $\{a_1, ..., a_n\}$. We say that $\phi$ *has remnant* if, for each $1 \leq i \leq n$, the word $a_i\phi$ has a nontrivial subword which does not cancel in any of the products

$$(a_j\phi)^{\pm 1}(a_i\phi), \quad (a_i\phi)(a_j\phi)^{\pm 1}, \tag{3.2}$$

except for $j = i$ with exponent $-1$. The maximal such noncancelling subword of $a_i\phi$ is called the remnant of $a_i$, denoted by $\text{Rem}_\phi(a_i)$.

The size of the remnant subwords can be measured in the following ways.

- If, for some natural number $l$, we have that $|\text{Rem}_\phi(a)| \geq l$ for all the generators $a$, then we say that $\phi$ has *remnant length* $l$.

- If, for some real number $r \in (0,1)$, it holds that $|\text{Rem}_\phi(a)| \geq r|a\phi|$ for all the generators $a$, then we say that $\phi$ has *remnant ratio* $r$.

**Example.**(see also [30]) The endomorphism $\phi$ of the free group on two generators $\{a, b\}$ given by the rule

$$\phi : \begin{cases} a & \longmapsto & a^2bab^{-2} \\ b & \longmapsto & ba^4ba^{-2} \end{cases}$$

has remnant. We have $\text{Rem}_\phi(a) = bab^{-1}$ and $\text{Rem}_\phi(b) = a^4b$.

Let us check that $\text{Rem}_\phi(a) = bab^{-1}$. It is easy to see that there is no cancellation in the products $(a\phi)(b\phi)^{-1}$, $(b\phi)^{-1}(a\phi)$, and $(a\phi)^2$. The longest terminal segment of $\phi(a)$ which cancels in

$$(a\phi)(b\phi) = (a^2bab^{-2})(ba^4ba^{-2})$$

is $b^{-1}$. Similarly, the word $a^2$ is the longest initial segment of $\phi(a)$ which cancels in

$$(b\phi)(a\phi) = (ba^4ba^{-2})(a^2bab^{-2}).$$

Therefore, $\text{Rem}_\phi(a) = bab^{-1}$. In the same way, one can show that $\text{Rem}_\phi(b) = a^4b$. Since $|\text{Rem}_\phi(g)| \geq 3$ for the generators $a$ and $b$, the remnant length of $\phi$ is 3. The remnant ratio of $\phi$ is $1/2$.

## 3.3 Membership in different doubly-twisted conjugacy classes

In this section we shall discuss two methods to detect whether arbitrary words belong to different doubly-twisted conjugacy classes with respect to homomorphisms $\phi$ and $\psi$. Both methods were developed by Staecker in [29] and [31]. The first one involves a joined remnant condition on $\phi$ and $\psi$, while the second one imposes a remnant condition on $\phi$ only so that the corresponding remnant words are longer than the images of the generators under $\psi$.

### 3.3.1 A strong remnant condition

Suppose that $\{a_1, ..., a_n\}$ is a generating set for the free group $G$. Denoting $G = \langle a_1, ..., a_n \rangle$, we consider the free product $G * G = \langle a_1, ..., a_n, a'_1, ..., a'_n \rangle$. There is a natural homomorphism $\phi * \psi : G * G \to H$ given by the rule

$$\phi * \psi : \begin{cases} a_i & \longmapsto & a_i\phi \\ a'_i & \longmapsto & a_i\psi. \end{cases}$$

Even though the equalizer subgroup is in general difficult to compute, the remnant property for homomorphisms $\phi, \psi$ forces it to be trivial. This is shown in the following lemma.

**Lemma 3.6.** *If $\phi * \psi : G * G \to H$ has remnant, then $(G\phi) \cap (G\psi) = \{1\}$. In particular, this means that $Eq(\phi, \psi) = \{1\}$.*

*Proof.* Suppose that $(G\phi) \cap (G\psi) \neq \{1\}$. Then there exist $x, y \in G$, both nontrivial, with $x\phi = y\psi$. So $(x\phi)(y\psi)^{-1} = 1$. However, this is impossible since $\phi * \psi$ has remnant; writing $x$ and $y$ in the generators shows that the word on the left hand side of the last equation cannot fully cancel. □

**Theorem 3.7** (Staecker, 2009). *Let $u$ and $v$ be different words in $H$. If $\phi^v * \psi$ has remnant, and if, for each generator $a$ of $G*G$, the remnant words $\text{Rem}_{\phi^v * \psi}(a)$ do not fully cancel in any product of the form*

$$(a\phi^v * \psi)v^{-1}u, \ u^{-1}v(a\phi^v * \psi), \qquad (3.3)$$

*then $[u] \neq [v]$.*

*Proof.* The idea of the proof is to show that the homomorphism $\widehat{\phi}^v * \widehat{\psi}$ (where $\widehat{\phi}$ and $\widehat{\psi}$ are defined as in Lemma 3.3) has remnant, and then to use Lemmas 3.3 and 3.6 to get the desired result. For a detailed treatment of the proof we refer the reader to [29]. □

It is also shown in [29] that if $\phi, \psi, u$ and $v$ are chosen at random, then $u$ and $v$ are not doubly-twisted conjugated with probability 1.

**Remark.** Notice that $\phi^v * id : G * G \to G$ can never have remnant. This follows easily from the fact that, when $v$ is written in the generators $a_i$, at least one $a_i$ $(= a_i' \phi^v * id)$ cancels in some product of the form (3.2). Therefore, Theorem 3.7 can never be adapted to distinguish twisted conjugacy classes.

### 3.3.2 A remnant inequality condition

Let us first introduce some notation. Let $A = \{a_1, ..., a_n\}$ be the set of generators of the free group $G$. Then any word $z \in G$ can be expressed in reduced form as

$$z = a_{j_1}^{\eta_1} \cdots a_{j_k}^{\eta_k}, \qquad (3.4)$$

where $a_{j_i} \in A$, $\eta_i = \pm 1$ for all $1 \leq i \leq k$, and every element of $A$ and its inverse are never adjacent.

Let $\phi : G \to H$ be a homomorphism between free groups, and $u, v \in H$. There is a natural homomorphism $\bar{\phi} = \phi * u * v : G * \mathbb{Z} * \mathbb{Z} \to H$ given by the rule

$$\phi * u * v : \begin{cases} a_i & \longmapsto & a_i\phi \\ b_1 & \longmapsto & u \\ b_2 & \longmapsto & v, \end{cases}$$

where $a_i$ $(1 \leq i \leq n)$ denote the generators of $G$, $b_1$ and $b_2$ the generators of the first and second $\mathbb{Z}$ factor, respectively.

**Theorem 3.8** (Staecker, 2009). *Let $\phi, \psi : G \to H$ be homomorphisms, and let $u, v \in H$, where the rank of $H$ is greater than 1. If $\bar{\phi}$ has remnant with*

$$|\text{Rem}_{\bar{\phi}}(a)| \geq |(a\psi)|$$

*for all the generators $a$ of $G$, then $[u] \neq [v]$.*

*Proof.* Assume that $[u] = [v]$. Then there is some $z \in G$ such that

$$z\psi = u^{-1}(z\phi)v. \tag{3.5}$$

We express $z$ as the reduced word $z = a_{j_1}^{\eta_1} \cdots a_{j_k}^{\eta_k}$, and write $X_i = (a_{j_i}\phi)^{\eta_i}$ and $Y_i = (a_{j_i}\psi)^{\eta_i}$. Then

$$u^{-1}(z\phi)v = u^{-1}X_1 \cdots X_k v = R_{u^{-1}}R_1 \cdots R_k R_v.$$

The right hand side of this equation is reduced. The words $R_{u^{-1}}$, $R_v$ and $R_i$ are the subwords of $u^{-1}$, $v$ and $X_i$ respectively which do not cancel in the above product. Since $\bar{\phi}$ has remnant, these words are nontrivial and $R_i$ contains $(\text{Rem}_{\bar{\phi}}(a_{j_i}))^{\eta_i}$ as a subword for all $1 \leq i \leq k$.

Comparing the length of $u^{-1}(z\phi)v$ and $z\psi$, we get a contradiction to (3.5). Indeed,

$$\begin{aligned}
|u^{-1}(z\phi)v| &= |R_u| + |R_v| + \sum_{i=1}^{k} |R_i| \\
&\geq |R_u| + |R_v| + \sum_{i=1}^{k} |\text{Rem}_{\bar{\phi}} a_{j_i}| \\
&> \sum_{i=1}^{k} |\text{Rem}_{\bar{\phi}} a_{j_i}| \\
&\geq \sum_{i=1}^{k} |Y_i| \\
&\geq |z\psi|.
\end{aligned}$$

$\square$

If $\psi$ is given and $\phi, u$ and $v$ are chosen at random, then the probability that $u$ and $v$ belong to different doubly-twisted conjugacy classes is 1. This is shown in [31].

Unlike Theorem 3.7, the previous theorem can be applied in the case where $\psi = id$.

**Corollary 3.9.** *Let $\phi : G \to G$ be an endomorphism of the free group $G$, and let $u, v \in G$. If $\phi * u * v$ has remnant, then $\text{TwOrb}(u) \neq \text{TwOrb}(v)$.*

*Proof.* The proof is an immediate consequence of Theorem 3.8. $\square$

## 3.4 Bounded solution length

We now present an algorithm that decides doubly-twisted conjugacy relations and provides a conjugating element if it exists [31]. This algorithm requires a remnant inequality condition similar to the one in Theorem 3.8. We also give a further condition for the equalizer subgroup to be trivial.

**Definition 3.10.** Given homomorphisms $\phi, \psi : G \to H$ and a pair $(u, v)$ of elements of $H$, we say that the pair $(u, v)$ has *bounded solution length* ($BSL$) if there is some $k > 0$ so that the equation $(z\phi)^{-1}u(z\psi) = v$ holds (if at all) only if $|z| \leq k$. The smallest such $k$ is called the *solution bound* ($SB$) for the pair $(u, v)$.

The next theorem states that any pair $(u, v)$ has $BSL$ with a predictable solution bound $SB$ whenever $\phi$ and $\psi$ satisfy a certain remnant inequality condition.

**Theorem 3.11** (Staecker, 2009). *Let $\phi, \psi : G \to H$ be homomorphisms such that $\phi$ has remnant with*

$$|\text{Rem}_\phi(a_i)| > |(a_i\psi)|$$

*for all the generators $a_i$ of $G = \langle a_1, ..., a_n \rangle$. Let $l := \min_{1 \leq i \leq n}(|\text{Rem}_\phi(a_i)| - |(a_i\psi)|)$. Any pair $(u, v)$ has $BSL$ with solution bound*

$$SB \leq \left\lfloor \frac{|u| + |v|}{l} \right\rfloor.$$

*Proof.* Let $u, v \in H$ and let $z \in G$ be a word of length $k$. The idea of the proof is to show that, for $k$ sufficiently large, we have that $(z\psi) \neq u^{-1}(z\phi)v$.

We express $z$ as the reduced word $z = a_{j_1}^{\eta_1} \cdots a_{j_k}^{\eta_k}$ (as described in Subsection 3.3.2), and write $X_i = (a_{j_i}\phi)^{\eta_i}$ and $Y_i = (a_{j_i}\psi)^{\eta_i}$. Then

$$u^{-1}(z\phi)v = u^{-1}X_1 \cdots X_k v = u^{-1}R_1 \cdots R_k v,$$

where each word $R_i$ is a subword of $X_i$ with $|R_i| \geq |\text{Rem}_\phi(a_{j_i})|$ (for all $1 \leq i \leq k$). The product $R_1 \cdots R_k$ is reduced but there might occur cancellation with $u^{-1}$ and $v$.

We now show that the words $u^{-1}(z\phi)v$ and $z\psi$ are of different length for sufficiently large $k$. Observe that

$$\begin{aligned}
|u^{-1}(z\phi)v| - |(z\psi)| &= |u^{-1}R_1\cdots R_k v| - |Y_1\cdots Y_k| \\
&\geq |R_1\cdots R_k| - |u| - |v| - |Y_1\cdots Y_k| \\
&= -|u| - |v| + (\sum_{i=1}^{k}|R_i|) - |Y_1\cdots Y_k| \\
&\geq -|u| - |v| + \sum_{i=1}^{k}(|R_i| - |Y_i|) \\
&\geq -|u| - |v| + \sum_{i=1}^{k}(|\text{Rem}_\phi a_{j_i}| - |Y_i|) \\
&\geq -|u| - |v| + kl.
\end{aligned}$$

The first inequality in the previous computation is an equality if both words $u^{-1}$ and $v$ fully cancel in the product $u^{-1}R_1\cdots R_k v$. The second inequality translates the fact that the product $Y_1\cdots Y_k$ might not be reduced. Finally, the last inequality follows by the assumption that $|\text{Rem}_\phi(a)| - |(a\psi)| \geq l$ for every generator $a$ of $G$.

Now, we may choose $k$ sufficiently large so that $|u^{-1}(z\phi)v| - |(z\psi)| > 0$; it suffices to choose $k > \frac{|u|+|v|}{l}$. For such a $k$ no solution can be found. This yields the desired solution bound

$$SB \leq \left\lfloor \frac{|u|+|v|}{l} \right\rfloor.$$

$\square$

The previous theorem yields the following algorithm for deciding doubly-twisted conjugacy relations. Suppose we want to compare the classes $[u]$ and $[v]$. If the hypotheses of Theorem 3.8 are satisfied, then we know that $[u] \neq [v]$. Otherwise, assuming that the remnant inequality condition of Theorem 3.11 holds, we have to verify whether $(z\psi) = u^{-1}(z\phi)v$ by a brute force check over all the words $z \in G$ with $|z| \leq \left\lfloor \frac{|u|+|v|}{l} \right\rfloor$.

Choosing the homomorphisms $\phi$ and $\psi$ at random, it is rather improbable that the condition of Theorem 3.11 holds. However, Staecker showed in [31] that this theorem can be applied whenever $\psi$ is fixed and $\phi$ chosen at random.

**Remark.** If we weaken the hypothesis of Theorem 3.11 to $|\text{Rem}_\phi(a_i)| = |(a_i\psi)|$ except for at least one of the generators (take $l \neq 0$), then we can at most deduce the distribution of the number of generators $a_i$ for which $|\text{Rem}_\phi(a_i)| \neq |(a_i\psi)|$ in a possible conjugating element $z$. In the proof of Theorem 3.11 we would get $k - m \leq \left\lfloor \frac{|u|+|v|}{l} \right\rfloor$ for words $z$ of length $k$ composed of $m$ generators for which $|\text{Rem}_\phi(a_i)| = |(a_i\psi)|$. In practice, this is rather useless (except maybe if we fix an upper bound for $|z| = k$).

**Example.** Let us consider the following pair of endomorphisms on the free group $F_2$ randomly generated by GAP [11] and stored for testing purpose

$$\phi : \begin{cases} a & \mapsto & b^2ab^2a^{-2} \\ b & \mapsto & a^{-1}b^{-1}ab^{-1}a^2b^{-1} \end{cases} \qquad \psi : \begin{cases} a & \mapsto & b^{-1} \\ b & \mapsto & ab^{-2}. \end{cases}$$

We can immediately make the following observations on $\phi$ and $\psi$

- the word $a\psi$ cancels in the product $(a\psi)(b\psi)^{-1}$, which implies that $\psi$ does not have remnant,

- the endomorphism $\phi$ has remnant with $\text{Rem}_\phi(a) = bab^2a^{-2}$ and $\text{Rem}_\phi(b) = a^{-1}b^{-1}ab^{-1}a^2$,

- $|\text{Rem}_\phi(g)| > |(g\psi)|$ for every generator $g \in \{a, b\}$ with $l = |\text{Rem}_\phi(b)| - |b\psi| = 3$.

So Theorem 3.7 does not allow to compare two classes of elements of $F_2$. But the hypotheses of Theorem 3.11 are indeed satisfied.

Recall that two elements $u, v \in F_2$ are said to be doubly-twisted conjugated if there exists an element $z \in F_2$ such that

$$(z\psi) = u^{-1}(z\phi)v. \tag{3.6}$$

We now analyze two selected examples in detail.

1. *Comparison of the classes* $[u] = \left[ba^{-1}ba^{-1}\right]$ *and* $[v] = \left[ba^{-1}b^{-2}\right]$. Notice that $\bar{\phi} := \phi * u * v$ does not have remnant, so that we cannot use Theorem 3.8. By Theorem 3.11, we know that $SB \leq \lfloor \frac{4+4}{3} \rfloor = 2$. A brute force check of all the words $z \in F_2$ of length at most 2 reveals that no such element verifies (3.6). Hence, we have that $[u] \neq [v]$.

2. *Comparison of the classes* $[u] = \left[a^{-1}b^{-1}\right]$ *and* $[v] = \left[ba^{-2}b^{-1}\right]$. Again Theorem 3.8 does not allow to make the comparison ($u$ cancels in the product $uv^{-1}$). So we have to check all the words $z \in F_2$ of length at most 2 since $SB \leq \lfloor \frac{2+4}{3} \rfloor = 2$ (by Theorem 3.11). This yields that $[u] = [v]$ with $b$ as conjugating element.

For the pair $\phi, \psi$ of endomorphisms, several tests to decide doubly-twisted conjugacy of two randomly chosen words $u$ and $v$ of $F_2$ (of a given maximal length) have been made. We refer the reader to Appendix A for the results of the tests as well as for the related GAP code.

The following corollary is an immediate consequence of Theorem 3.11.

**Corollary 3.12.** *Suppose that $G$ and $H$ are finitely generated free groups. If $\phi : G \to H$ is a homomorphism with $|\text{Rem}_\phi(a)| > 1$ for every generator $a \in G$, then there exists an algorithm to decide, for any $u, v \in H$, whether or not*

$$v = (z\phi)^{-1}uz$$

*for some $z \in G$ ($\phi$-twisted conjugacy).*

Our next proposition gives another necessary condition for the equalizer subgroup to be trivial.

**Proposition 3.13.** *Let $\phi, \psi : G \to H$ be homomorphisms. Suppose that $\phi$ has remnant with*
$$|\text{Rem}_\phi(a_i)| > |(a_i\psi)| \tag{3.7}$$
*for all the generators $a_i$ of $G = \langle a_1, ..., a_n \rangle$, except for one of them for which equality in (3.7) holds. If in addition the images of this generator under $\phi$ and $\psi$ are different, then $Eq(\phi, \psi) = \{1\}$.*

*Proof.* Let $l := \min_{1 \leq i \leq n} (|\text{Rem}_\phi(a)| - |(a\psi)|) \neq 0$, and let $a_r$ $(1 \leq r \leq n)$ be the generator for which $|\text{Rem}_\phi(a_r)| = |(a_r\psi)|$. Let $z \neq 1$ be an arbitrary word in $G$ of length $k$ which contains $m$ $(0 \leq m \leq k)$ occurrences of $a_r^{\pm 1}$.

Expressing $z$ as the reduced word $z = a_{j_1}^{\eta_1} \cdots a_{j_k}^{\eta_k}$, we write $X_i = (a_{j_i}\phi)^{\eta_i}$ and $Y_i = (a_{j_i}\psi)^{\eta_i}$. Then
$$z\phi = X_1 \cdots X_k = R_1 \cdots R_k,$$
where each word $R_i$ is a subword of $X_i$ with $|R_i| \geq |\text{Rem}_\phi(a_{j_i})|$ (for all $1 \leq i \leq k$).

We have
$$\begin{aligned}
|(z\phi)| - |(z\psi)| &= |R_1 \cdots R_k| - |Y_1 \cdots Y_k| \\
&\geq \sum_{i=1}^k |R_i| - \sum_{i=1}^k |Y_i| \\
&\geq \sum_{i=1}^k (|\text{Rem}_\phi(a_{j_i})| - |Y_i|) \\
&\geq (k-m)l.
\end{aligned}$$

Now $|(z\phi)| - |(z\psi)| > 0$ if and only if $k > m$. So if the element $z$ is in $Eq(\phi, \psi)$, then $z = a_r^n$ for some $n \in \mathbb{Z}$. Since $(a_r\phi)^n = (a_r\psi)^n$ implies that $a_r\phi = a_r\psi$ (which does not hold by our assumption), we conclude that $Eq(\phi, \psi) = \{1\}$. □

We obtain the following corollary as an immediate consequence of Proposition 3.13.

**Corollary 3.14.** *Let $\phi, \psi : G \to H$ be homomorphisms. If $\phi$ has remnant with*
$$|\text{Rem}_\phi(a_i)| > |(a_i\psi)|$$
*for all the generators $a_i$ of $G = \langle a_1, ..., a_n \rangle$, then $Eq(\phi, \psi) = \{1\}$.*

*Proof.* Let $l := \min_{1 \leq i \leq n} (|\text{Rem}_\phi(a)| - |(a\psi)|)$. It holds for all words $z \neq 1$ in $G$ that
$$|(z\phi)| - |(z\psi)| \geq kl > 0,$$
where $k$ denotes the length of $z$. Consequently, $Eq(\phi, \psi) = \{1\}$. (For more details see proof of Proposition 3.13) $\square$

**Remark.** Consider the case where $|\text{Rem}_\phi(a_r)| - |(a_r\psi)| = d < 0$ for exactly one generator $a_r$ of the free group $G$ (for the other generators the difference is supposed to be strictly greater than 0). If there exists an element $z \in Eq(\phi, \psi)$ of length $k$, then the number of occurrences of $a_r^{\pm 1}$ in $z$ must be $\frac{kl}{l-d}$ ($\geq 1$).

## 3.5 Open problems

In this final section we address some interesting problems related to the equalizer subgroup which, to our knowledge, are still open at the time of writing. We also give a comment on some of them. In the following we denote by $\phi, \psi : G \to H$ homomorphisms between free groups.

**Question 1.** Is there a relation between the set $C(\phi, \psi) := \{w \in G \mid [w\phi, w\psi] = 1\}$ and the subgroup $Eq(\phi, \psi)$?

Clearly, $Eq(\phi, \psi) \subseteq C(\phi, \psi)$. Conversely, suppose that $\tilde{w} \in C(\phi, \psi)$. Then $\langle \tilde{w}\phi, \tilde{w}\psi \rangle = \langle u \rangle$ for some computable $u \in H$. So $\tilde{w} \in Eq(\phi, \psi)$ if only if $\tilde{w}\phi = u = \tilde{w}\psi$.

**Question 2.** Is there an algorithm for deciding whether $Eq(\phi, \psi)$ is trivial or not?

As already mentioned this problem is equivalent to PCP in free groups. Several decidability and undecidability results on this problem are known in the case where $\phi, \psi : A^* \to B^*$ are homomorphisms between free monoids. For instance, if $n := |A| = 2$ (number of letters in the alphabet $A$), then PCP is decidable [9]. However, if $n \geq 7$, then the problem is undecidable [23]. It is not yet known whether or not it is decidable for $3 \leq n \leq 6$. By restricting $\phi, \psi$ to be injective, PCP remains undecidable [21]. What results can be obtained in the world of free groups?

**Question 3.** Given a finitely generated subgroup $S$ of the free group $G$, do there exist homomorphisms $\phi, \psi$ such that $Eq(\phi, \psi) = S$?

**Question 4.** Given $w \in Im(\phi) \cap Im(\psi)$, can we decide whether or not there exists $z \in Eq(\phi, \psi)$ such that $z\phi = z\psi = w$?

In the case where $\phi, \psi$ are automorphisms, this is trivially decidable. It suffices to compare $w\phi^{-1}$ to $w\psi^{-1}$.

Another particular situation occurs when $\phi$ has remnant (which implies that $\phi$ is injective). Suppose $w \in Im(\phi) \cap Im(\psi)$. Then there exist words $z, \tilde{z} \in G$ such that $z\phi = w$ and $\tilde{z}\psi = w$. Since $\phi$ has remnant, we have that $|(z\phi)| - |w| \geq |z| - |w|$. So the length of $z$ is at most $|w|$. Now $z \in Eq(\phi, \psi)$ if and only if $z$ belongs to the coset $\tilde{z}K$, where $K$ denotes the kernel of $\psi$. However there is no possibility to get any information about the coset representative $\tilde{z}$ so that a brute force check over the elements in $\tilde{z}K$ of length at most $|w|$ is impossible.

**Question 5.** Can we determine all endomorphisms $\phi, \psi$ of $F_2$ with $Fix(\phi) \cap Fix(\psi) = Eq(\phi, \psi)$? Can we restate the same question by replacing $Fix(\phi) \cap Fix(\psi)$ with the commutator subgroup $\langle [G\phi, G\psi] \rangle$?

Using the fact that $Fix(\phi) \cap Fix(\psi) \leq Eq(\phi, \psi) \leq \langle [G\phi, G\psi] \rangle$, these questions might be justified. Below we analyse a few particular cases.

Suppose that $\phi, \psi$ are non-trivial endomorphisms of $F_2 = \langle a, b \rangle$ with $rk(Fix(\phi) \cap Fix(\psi)) = 2$. Then, according to Theorem 3.4 of [35], it holds that $\phi, \psi \in Aut(F_2)$ and $Fix(\phi) \cap Fix(\psi) = \langle a, bab^{-1} \rangle = Fix(\phi) = Fix(\psi)$. Also $Eq(\phi, \psi) = Fix(\psi^{-1}\phi) = Fix(\phi^{-1}\psi)$. Do such endomorphisms exist?

Let $\phi$ be an endomorphism of $F_2$. If $\phi$ is not bijective, then $rk(Fix(\phi)) \leq 1$ [33]. By a result of Hanna Neumann (see, for example, [3]), we deduce that in this case $rk(Fix(\phi) \cap Fix(\psi)) \leq 2(rk(Fix(\phi)) - 1)(rk(Fix(\psi)) - 1) + 1 \leq 1$. So for all pairs of endomorphisms, with one of them being non-bijective, that satisfy the above equality, the equalizer subgroup is either trivial or infinite cyclic.

# Chapter 4

# Application: An Authentication Scheme

The main goal of an authentication protocol is to allow a legitimate user to prove his identity over an insecure channel to a server using his private key without leaking any information about the key. We discuss in this chapter an authentication scheme introduced by Shpilrain and Ushakov in [28] that is based on the *doubly-twisted conjugacy search problem*:

> Given a pair of endomorphisms $\phi, \psi$ of a group $G$ and a pair of elements $w, t \in G$, find an element $s \in G$ such that $t = (s^{-1}\phi)w(s\psi)$.

The platform $G$ will be the semigroup of all $2 \times 2$ matrices over truncated one-variable polynomials over $\mathbb{F}_2$, the finite field of two elements. In this platform semigroup we consider matrix transposition instead of inversion, so that the above equality becomes $t = (s^T\phi)w(s\psi)$. To conclude the chapter, we present a heuristic attack on this authentication scheme that has been recently suggested by Grassl and Steinwandt in [14]. This attack relies on efficient Gröbner bases computations and on simple elimination techniques for linear equations.

## 4.1 The protocol

The group-based authentication protocol that we describe in this section is due to Shpilrain and Ushakov [28]. As the Fiat-Shamir protocol (see, for instance, [34]), this protocol involves repeating several times a three-pass challenge-response step. Here Alice is referred to as the prover and Bob is referred to as the verifier.

*Setup of public parameters:* Let $G$ be a platform semigroup and $*$ be an antihomomorphism of $G$ into itself (that is, $(ab)^* = b^*a^*$ for any $a, b \in G$).

Assume that the doubly-twisted conjugacy search problem in $G$ is hard. We let $k$ be a security parameter.

*Setup of the keys:* Alice chooses a secret key $s \in G$. She publishes a pair of endomorphisms $\phi, \psi$ of $G$ and two elements $w, t \in G$ such that $t = (s^*\phi)w(s\psi)$.

*Protocol:* Perform $t$ rounds in which

1. Alice picks a random[1] $r \in G$ and sends a commitment $u = (r^*\phi)t(r\psi)$ to Bob.

2. Bob chooses a random bit $c$ and sends it as a challenge to Alice.

3. If $c = 0$, then Alice sends $v = r$ to Bob. Bob aborts the protocol and rejects unless the equality $u = (v^*\phi)t(v\psi)$ is satisfied.

4. If $c = 1$, then Alice sends $v = sr$ to Bob. Bob aborts the protocol and rejects unless the equality $u = (v^*\phi)w(v\psi)$ is satisfied.

After $k$ successful rounds, Bob accepts the authentication.

If Alice and Bob behave as specified, then the final check of each round leads

- in Step 3. to $(v^*\phi)t(v\psi) = (r^*\phi)t(r\psi) = u$ since $v = r$, and

- in Step 4. to $(v^*\phi)w(v\psi) = ((sr)^*\phi)w((sr)\psi) = ((r^*s^*)\phi)w((sr)\psi) = (r^*\phi)(s^*\phi)w(s\psi)(r\psi) = (r^*\phi)t(r\psi) = u$ since $v = sr$.

So a correct answer of the prover leads to acceptance by the verifier.

An adversary Eve can impersonate Alice with probability $\frac{1}{2^k}$ if she behaves as follows. In each round, Eve predicts the challenge $\tilde{c}$ and picks a random $r \in G$. If $\tilde{c} = 0$, then she sends $u = (r^*\phi)t(r\psi)$ to Bob. If $\tilde{c} = 1$, then she sends $u = (r^*\phi)w(r\psi)$ to Bob. If the challenge $c$ chosen by Bob is equal to $\tilde{c}$, then Eve can answer it by $v = r$. Otherwise it fails.
Alternatively, Eve could try to find a secret key $\tilde{s}$ (not necessarily equal to the secret key $s$) from the knowledge of $w$, $t$, $\phi$ and $\psi$ such that $t = (\tilde{s}^*\phi)w(\tilde{s}\psi)$, that is, she has to find a solution to the doubly-twisted conjugacy search problem.

## 4.2 The platform

The suggested platform $G$ for the protocol is the semigroup of all $2 \times 2$ matrices over truncated one-variable polynomials over $\mathbb{F}_2$, the field of two elements. Truncated, more precisely $N$-truncated, one-variable polynomials over $\mathbb{F}_2$ are expressions of the form $\sum_{0 \leq i \leq N-1} a_i x^i$, where $a_i$ are elements of $\mathbb{F}_2$ and $x$ is a variable. In other words, $N$-truncated polynomials are elements of $R := \mathbb{F}_2[x]/\langle x^N \rangle$, the factor algebra of the algebra $\mathbb{F}_2[x]$ of one-variable polynomials over $\mathbb{F}_2$ by the ideal generated by $x^N$. This set can be identified by the set of all polynomials of degree strictly smaller than $N$. Moreover, we denote by $R^* := \{f(x) \in R \mid f(0) \neq 0\}$ the set of all polynomials in $R$ with constant term 1. For polynomials in $R$, we have the following operations:

$$\begin{array}{lll} \text{Addition:} & (f(x), g(x)) \mapsto & f(x) + g(x) \\ \text{Multiplication:} & (f(x), g(x)) \mapsto & f(x)g(x) \bmod x^N \\ \text{Composition:} & (f(x), g(x)) \mapsto & f \circ g := f(g(x)) \bmod x^N. \end{array}$$

---
[1] All random choices in this chapter are with respect to uniform distribution.

Endomorphisms of $G$ are naturally induced by endomorphisms of the algebra of truncated polynomials $R$. Any map of the form $x \mapsto p(x) \in R \setminus R^*$ (that is, $p(x)$ has zero constant term) can be extended to an endomorphism $\phi_p$ of $R$ by the rule $\phi_p(f(x)) = f(p(x)) \bmod x^N$, where $f \in R$. The condition that $p$ has zero constant term ensures that the ideal generated by $x^N$ is invariant under $\phi_p$, or equivalently, that $\phi_p$ is indeed an endomorphism of $R$. To check this, consider the polynomial $p(x) = \sum_{i=1}^{N-1} a_i x^i$. We have that

$$\begin{aligned} \phi_p(x^2) &= (p(x))^2 \bmod x^N \\ &= (\sum_{i=1}^{N-1} a_i x^i)(\sum_{j=1}^{N-1} a_j x^j) \bmod x^N \\ &= \sum_{i=2}^{(N-1)^2} (\sum_{k=1}^{i-1} a_k a_{i-k}) x^i \bmod x^N \\ &= \begin{cases} \sum_{i=1}^{(N-1)/2} a_i x^{2i} & \text{if } N-1 \text{ even} \\ \sum_{i=1}^{(N-2)/2} a_i x^{2i} & \text{if } N-1 \text{ odd} \end{cases} \\ &\in \langle x^2 \rangle \end{aligned}$$

Clearly, $\phi_p(x^N) = (p(x))^N \in \langle x^N \rangle$ by a similar computation as above.

Using elementary operations on polynomial functions, it can be easily shown that $\phi_p$ is both an additive and a multiplicative endomorphism of $R$. Then, it extends naturally to an endomorphism of the semigroup $G$ of all $2 \times 2$ matrices over truncated one-variable polynomials.

## 4.3 The parameters and key generation

The parameters and key generation of the protocol can be summarized as follows.

*Setup of public parameters:* Set $N = 300$. The platform $G$ is the semigroup of all $2 \times 2$ matrices over $R := \mathbb{F}_2[x]/\langle x^{300} \rangle$. The antihomomorphism of $G$ is given by matrix transposition $((AB)^T = B^T A^T$ for any $A, B \in G)$.

*Setup of the keys:* Alice's secret key consists of 4 random polynomials $s_i \in R^*$ for $i = 1, ..., 4$. Her public key consists of polynomials $p, q \in R \setminus R^*$ and $w_i \in R^*$, for $i = 1, ..., 4$, chosen uniformly at random, and $t_i \in R$ ($i = 1, ..., 4$) where

$$\begin{pmatrix} t_1 & t_2 \\ t_3 & t_4 \end{pmatrix} := \begin{pmatrix} s_1 \circ p & s_3 \circ p \\ s_2 \circ p & s_4 \circ p \end{pmatrix} \begin{pmatrix} w_1 & w_2 \\ w_3 & w_4 \end{pmatrix} \begin{pmatrix} s_1 \circ q & s_2 \circ q \\ s_3 \circ q & s_4 \circ q \end{pmatrix}$$

We illustrate the above protocol by an example and refer the reader to Appendix B for the related GAP code.

**Example.** In our example polynomials are truncated at $N = 7$ so that we consider the ring $R = \mathbb{F}_2[x]/\langle x^7 \rangle$. Note that in the GAP code the secret key $s$ as well as the randomly chosen $r$ are given by matrices of polynomial functions.

*Secret key:*

$$s = \begin{pmatrix} x^4 + x^3 + x^2 + x + 1 & x^6 + x^4 + x^2 + 1 \\ x^6 + x + 1 & x^2 + 1 \end{pmatrix}$$

*Public key:*

$$\phi_p(x) = p(x) = x^6 + x^4 + x^3 + x^2$$
$$\psi_q(x) = q(x) = x^3 + x^2 + x$$
$$w = \begin{pmatrix} x^3 + x^2 + 1 & x^3 + x + 1 \\ x^5 + x^3 + 1 & x^6 + x^5 + x^4 + x^3 + x^2 + 1 \end{pmatrix}$$
$$t = \begin{pmatrix} x^6 + x^4 + x^2 + x & x \\ x^5 + x^4 + x^3 + x^2 + x & x^6 + x^4 + x^3 + x \end{pmatrix}$$

*One round of the protocol:*

1. Alice picks randomly

$$r = \begin{pmatrix} x^5 + x^4 + x^3 + x^2 + 1 & x^6 + x^5 + x^4 \\ x^2 & x^3 + x^2 + x \end{pmatrix}$$

   and sends the commitment $u = (r^T \circ \phi_p) t (r \circ \psi_q)$ to Bob, that is,

$$u = \begin{pmatrix} x^5 + x^4 + x^2 + x & x^5 + x^2 \\ x^3 & x^6 + x^5 + x^4 \end{pmatrix}.$$

2. Bob chooses a random bit $c$ and sends it as a challenge to Alice.

3. If $c = 0$, then Alice sends $v = r$ to Bob. Bob computes

$$(v^T \circ \phi_p) t (v \circ \psi_q) = \begin{pmatrix} x^5 + x^4 + x^2 + x & x^5 + x^2 \\ x^3 & x^6 + x^5 + x^4 \end{pmatrix}$$

   compares it to $u$ and accepts!

4. If $c = 1$, then Alice sends $v = sr$ to Bob, that is,

$$v = \begin{pmatrix} x^6 + x^4 + x^3 + x^2 + x + 1 & x^2 + x \\ x^4 + x + 1 & x^5 + x^2 + x \end{pmatrix}.$$

Bob computes

$$(v^T \circ \phi_p) w (v \circ \psi_q) = \begin{pmatrix} x^5 + x^4 + x^2 + x & x^5 + x^2 \\ x^3 & x^6 + x^5 + x^4 \end{pmatrix}$$

compares it to $u$ and accepts!

## 4.4 Cryptanalysis

In this section we present a heuristic attack on the previously introduced authentication scheme that is due to Grassl and Steinwandt [14]. A review on cryptanalysis using Gröbner bases is given in [10]. Some background information on Gröbner bases and monomial orderings is taken from this paper.

Recall that the public key consists of polynomials $p, q \in R \setminus R^*$, $w_i \in R^*$, and $t_i \in R$ ($1 \leq i \leq 4$). The goal of the attack is to detect polynomials $\tilde{s}_i \in R^*$ (not necessarily equal to the actual secret key $s_i$, for $i = 1, ..., 4$) with

$$\begin{pmatrix} t_1 & t_2 \\ t_3 & t_4 \end{pmatrix} = \begin{pmatrix} \tilde{s}_1 \circ p & \tilde{s}_3 \circ p \\ \tilde{s}_2 \circ p & \tilde{s}_4 \circ p \end{pmatrix} \begin{pmatrix} w_1 & w_2 \\ w_3 & w_4 \end{pmatrix} \begin{pmatrix} \tilde{s}_1 \circ q & \tilde{s}_2 \circ q \\ \tilde{s}_3 \circ q & \tilde{s}_4 \circ q \end{pmatrix} \quad (4.1)$$

To begin with, we write

$$\tilde{s}_i = 1 + \sum_{j=1}^{N-1} y_{ij} x^j \ (1 \leq i \leq 4, \ 1 \leq j < N),$$

where the variables $y_{ij}$ are indeterminates over $\mathbb{F}_2$. A system of $4(N-1)$ polynomial equations in $4(N-1)$ unknowns is obtained in the following way. First, we evaluate the matrix product on the right hand-side of (4.1) with each $\tilde{s}_i$ replaced by the sum $1 + \sum_{j=1}^{N-1} y_{ij} x^j$ and computations performed modulo $x^N$. Then, for each matrix entry, we equate the coefficients of all the $x^j$ ($1 \leq j < N$) on both sides of (4.1).

To the resulting system of polynomial equations under the form

$$f_1(y_{11}, ..., y_{4(N-1)}) = ... = f_{4(N-1)}(y_{11}, ..., y_{4(N-1)}) = 0$$

we associate the ideal generated by the set of polynomials $P = \{f_1, ..., f_{4(N-1)}\}$. To this set $P$ we add the *field polynomials* $y_{ij}^2 - y_{ij}$ ($1 \leq i \leq 4$, $1 \leq j < N$). This ensures that the ideal generated by $P$ is 0-dimensional (finite number of solutions).

Gröbner bases provide an important tool for solving polynomial systems of equations. To avoid the cost of a full Gröbner basis computation, the computations in the attack involve truncated degree-2 Gröbner bases obtained by disregarding $S$-polynomial pairs whose total degree is greater than 2. Astonishingly these truncated Gröbner bases allow the derivation of linear equations that suffice to attack the authentication scheme in question.

The monomial ordering used for the truncated Gröbner bases computations is the graded reverse lex (grevlex) ordering. Recall that the grevlex ordering is defined as follows. Let the variables in a monomial $y$ be ordered by $y_{11} > y_{12} > ... > y_{4(N-1)}$. For two monomials $y^\alpha$ and $y^\beta$ (with $\alpha, \beta \in \mathbb{Z}_{\geq 0}^{4(N-1)}$), $y^\alpha >_{grevlex} y^\beta$ if $deg(y^\alpha) > deg(y^\beta)$ or $deg(y^\alpha) = deg(y^\beta)$ and, in $\alpha - \beta \in \mathbb{Z}^{4(N-1)}$, the right-most nonzero entry is negative. This ordering is less intuitive than the well-known lexicographic ordering, but it is proven to be generally the most efficient ordering for computing Gröbner bases (in particular in terms of time complexity).

The heuristic attack on the secret key proceeds as follows.

1. Initialize $B := \left\{ y_{ij}^2 - y_{ij},\ 1 \leq i \leq 4,\ 1 \leq j < N \right\}$ and $d := 1$.

2. Equate the coefficients of $x^d$ on both sides of (4.1). This leads to a set of polynomials $B_d$. Set $B := B \cup B_d$.

3. Compute the truncated degree-2 Gröbner basis $G$ of $B$ with respect to the grevlex ordering.

4. Extract all linear polynomials (that is, polynomials of total degree 1) from $G$ and compute a reduced row echelon form.

   - From polynomials of the form $y_{i_0 j_0}$ or $y_{i_0 j_0} - 1$ we can uniquely determine some coefficients of the (alternative) secret key.
   - The remaining linear polynomials are of the form
     $y_{i_0 j_0} - \sum_{(i,j) \neq (i_0, j_0)} y_{ij}$ with $y_{i_0 j_0}$ not occurring in the other linear polynomials. We reduce the number of indeterminates by replacing each occurrence of $y_{i_0 j_0}$ in $B$ with the corresponding sum $\sum_{(i,j) \neq (i_0, j_0)} y_{ij}$.

5. If $d < N - 1$ and the complete (alternative) secret key has not yet been discovered, then we increase $d$ by 1 and return to Step 1.

6. The variables in the remaining polynomial equations are set to 0 and we check whether or not the resulting candidate key verifies (4.1).

The heuristic is based on the observation that in many cases the only constraints on the remaining variables are of the form $y_{ij}^2 = y_{ij}$, so that one can choose any value from $\mathbb{F}_2$ for these $y_{ij}$.

For the proposed parameter value $N = 300$, experimental results in the Magma language[2] have led the authors of [14] to a success rate of 88% within a few minutes with moderate resources. Hence the above authentication scheme does not provide strong cryptographic security guarantees.

---

[2] W. Bosma, J. Cannon, and C. Playoust. The Magma algebra system. I. The user language. *J. Symbolic Comput.* **24** (3-4) (1997), 235-265.

# Appendix A

# Implementation in GAP

The goal of this appendix is to present the implementation and testing of functions that appear in Staecker's Theorems 3.7 and 3.11. To this end, some auxiliary programs from the folder ntwm2.tar (current release: v2.0) on Staecker's homepage[1] are required.

The GAP code is organized in two files. The first file "rembsl.gap" contains programs from the files "maps.gap", "brute.gap", "bsl.gap", "wagner.gap" of the folder ntwm2.tar. Relevant modifications are indicated by dashed lines. The second file "dtc.gap" includes the code for testing the hypotheses of Theorems 3.7 and 3.11. Moreover the algorithm described in Section 3.4 is implemented here.

We first describe the functions of both files in Section A.1, then list the code in Section A.2, and finally generate some tests on the command line in Section A.3.

## A.1 Documentation

This documentation should give the interested reader the necessary information to predict the outcome of each of the defined functions.

The terminology used to characterize the input to these functions is as follows:

- $G, H$: finitely generated free groups
- $GFF$: free product $G * \mathbb{Z} * \mathbb{Z}$ (as defined in (3.3.2))
- $f, g$: homomorphisms between finitely generated free groups
- $x$: generator of some finitely generated free group
- $u, v, w$: freely reduced words

---

[1] http://192.160.243.156/~cstaecker/research.html (Accessed 09/10/2009)

- $k, l, m, n$: positive integers.

The file "rembsl.gap" includes the following functions.

rf( $G$, $H$, $l$) – Generates a random homomorphism between $G$ and $H$ with word lengths at most $l$.

InitialCancellingSegment( $w$, $u$, $G$) – Returns the initial segment of $w$ that cancels in the product $uw$.

TerminalCancellingSegment( $w$, $u$, $G$) – Returns the terminal segment of $w$ that cancels in the product $wu$.

RemnantTriple( $f$, $x$, $G$, $H$) – Returns a triple $[p, r, s]$, where $p$ is the longest initial segment of $xf$ that cancels in the product $(yf)^{\pm 1}(xf)$, $s$ the longest terminal segment of $xf$ that cancels in the product $(xf)(yf)^{\pm 1}$ and $r$ the possible remainder (remnant) (excluded: $y = x$ with exponent $-1$). If $|p| + |s| \geq |(xf)|$ (indicating that $f$ has no remnant), then the triple $[xf, One(G), One(G)]$ is returned.

RemnantOfGenerator( $f$, $x$, $G$, $H$) – Returns the remnant of the generator $x$ under the homomorphism $f$.

HasRemnant( $f$, $G$, $H$) – Indicates whether or not the homomorphism $f$ has remnant. If, for at least one of the generators of $G$, the function RemnantOfGenerator returns $One(G)$, then HasRemnant returns false, otherwise it returns true.

HasBSL( $f$, $g$, $G$, $H$) – Indicates whether or not the $BSL$ condition on $f$ and $g$ holds. Recall that the $BSL$ condition is satisfied whenever $|\text{Rem}_f(a)| > |(ag)|$ for all the generators $a$ of $G$.

DoubleTwistConj( $f$, $g$, $v$, $w$) – Computes the doubly-twisted conjugacy relation on $v$, that is, $(w^{-1}f)v(wg)$.

IsDoublyTwistedConjugateBruteForce( $f$, $g$, $u$, $v$, $G$, $H$, $k$) – Checks iteratively whether or not the relation $v = (wf)^{-1}u(wg)$ holds for some $w \in G$ of length at most $k$. If no such element can be found, the function returns fail. Note that in the worst case $1 + \sum_{i=1}^{k} 2m(2m-1)^{i-1}$ checks of the doubly-twisted conjugacy relation are performed ($m$ denotes the rank of the free group $G$), so that for large values of $k$ the computation risks to be very slow.

In the file "dtc.gap" we define some free groups, homomorphisms and words as global variables (mainly for testing purpose). Note that the free groups $F11$ and $F12$ are used in some programs to define the free product $G * \mathbb{Z} * \mathbb{Z}$. Further, the file contains the following functions.

rel( $G$, $l$) – Generates a random element of the free group $G$ of length at most $l$.

FreeProductHomomorphism( $f$, $u$, $v$, $G$, $H$, $GFF$) – Let $GFF$ be the free product $G * F11 * F12$. FreeProductHomomorphism computes the homomorphism $\bar{f} = f * u * v : GFF \to H$ as defined in (3.3.2).

HasBSLW( $f$, $g$, $u$, $v$, $G$, $H$) – Indicates whether or not the condition of Theorem 3.8 holds. Recall that this condition is satisfied whenever $f * u * v$ has remnant with $|\text{Rem}_{\bar{f}}(a)| \geq |(ag)|$ for all the generators $a$ of $G$.

MinBSL( $f$, $g$, $G$, $H$) – If the $BSL$ condition holds, then the function MinBSL returns the minimum of the difference $|\text{Rem}_f(a)| - |(ag)|$ ($> 0$) over all the generators $a$ of $G$. Otherwise MinBSL returns fail.

IsDTCBSL( $f$, $g$, $u$, $v$, $G$, $H$) – If the $BSL$ condition is verified, then IsDTCBSL provides an algorithm for deciding doubly-twisted conjugacy of the words $u$ and $v$. Before we apply the procedure of Theorem 3.11, we check equality of the words $u, v$ and also HasBSLW (in case one of these tests returns true we already know the outcome). Otherwise an error message is returned.

DTCBSLTests( $f$, $g$, $G$, $H$, $m$, $n$) – In case the $BSL$ condition is verified, DTCBSLTests generates random pairs of distinct elements from the free group $H$ of lengths at most $m$ which satisfy the doubly-twisted conjugacy relation for some $w \in G$ (out of $n$ trials).

## A.2 Code listing

### A.2.1 "rembsl.gap"

```
# File name: rembsl.gap
# Author: Christopher Staecker
# Last modified: 24/05/2009
# Remark: Relevant code modifications by the author of this
#         thesis are indicated between dashed lines

#*************************Terminology*************************
# G,H: finitely generated free groups
# f,g: homomorphisms between finitely generated free groups
# x,u,v,w: freely reduced words
# k,l: positive integers
#*************************************************************
```

```
# rf( G, H, l) generates a random homomorphism between the free
# groups G and H with word lengths at most l
rf := function( G, H, l)
  local letterlist, ims, gens;
  gens := GeneratorsOfGroup(H);
  letterlist := function()
        return List([1..l], n -> Random(gens)^Random([-1,1]));
  end;
  ims := List(GeneratorsOfGroup(G),g -> Product(letterlist()));
  return GroupHomomorphismByImagesNC(G,H,GeneratorsOfGroup(G),
                                     ims);
end;

# InitialCancellingSegment( w, u, G) returns the initial seg-
# ment of w that cancels in the product uw
InitialCancellingSegment := function( w, u, G)
  local a,b;
  if (w = One(G) or u = One(G)) then
    return One(G);
  fi;
  a := Subword(w,1,1);
  b := Subword(u^-1, 1, 1)^-1;
  if (a*b = One(G)) then
    return a * InitialCancellingSegment(a^-1 * w, u * b^-1, G);
  else
    return One(G);
  fi;
end;

# TerminalCancellingSegment( w, u, G) returns the terminal seg-
# ment of w that cancels in the product wu
TerminalCancellingSegment := function(w, u, G)
  return (InitialCancellingSegment(w^-1, u^-1, G))^-1;
end;

# RemnantTriple( f, x, G, H) returns a triple [p,r,s], where p
# and s are the subwords of (xf) respectively before and after
# the remnant subword (if it exists). If |p|+|s| >= |(xf)|
# (indicating that f has no remnant), then the triple
# [xf,One(G),One(G)] is returned.
```

```
RemnantTriple := function( f, x, G, H)
  local gpm, p, s, r;
  gpm := Concatenation(GeneratorsOfGroup(G),
    List(GeneratorsOfGroup(G),g -> g^-1));
  gpm := Filtered(gpm,g -> g <> x^-1);
  p := Maximum(List(gpm,g -> InitialCancellingSegment(x^f,g^f,
                                                     H)));
  s := Maximum(List(gpm,g -> TerminalCancellingSegment(x^f,g^f,
                                                      H)));
#----------------------------------------------------------------
  if (Length(p)+Length(s) >= Length(x^f)) then
    return [x^f,One(G),One(G)];
#----------------------------------------------------------------
  else
    r := p^-1 * x^f * s^-1;
    return [p,r,s];
  fi;
end;

# RemnantOfGenerator( f, x, G, H) returns the remnant of the
# generator x under the homomorphism f : G -> H
RemnantOfGenerator := function( f, x, G, H)
  return RemnantTriple(f,x,G,H)[2];
end;

# HasRemnant( f, G, H) indicates whether or not the homomor-
# phism f : G -> H has remnant
HasRemnant := function( f, G, H)
  local a;
  for a in GeneratorsOfGroup(G) do
    if Length(RemnantOfGenerator(f,a,G,H)) = 0 then
      return false;
    fi;
  od;
  return true;
end;

# HasBSL( f, g, G, H) indicates whether or not
# |Rem_{f}(a)| > |(ag)| for all the generators a of G
HasBSL := function( f, g, G, H)
  local a, r, l;
```

```
#----------------------------------------------------------------
  if not HasRemnant(f,G,H) then
  return false;
#----------------------------------------------------------------
  else
  for a in GeneratorsOfGroup(G) do
    r := RemnantOfGenerator(f,a,G,H);
    l := a^g;
    if (Length(r) <= Length(l)) then
      return false;
      fi;
    od;
    return true;
  fi;
end;

# DoubleTwistConj( f, g, v, w) computes the doubly-twisted
# conjugacy relation on v, that is, (w^{-1}f)v(wg)
DoubleTwistConj := function( f, g, v, w)
#----------------------------------------------------------------
  return (w^-1)^f * v * w^g;
#----------------------------------------------------------------
end;

# IsDoublyTwistedConjugateBruteForce( f, g, u, v, G, H, k)
# checks iteratively whether or not the relation
# v=(wf)^{-1}u(wg) holds for some w in G of length at most k
IsDoublyTwistedConjugateBruteForce := function( f, g, u, v, G,
H, k)
  local I, w;
  I := Iterator(G);
  w := NextIterator(I);
  while Length(w) <= k do
    if DoubleTwistConj(f,g,u,w) = v then
      Print("The elements ",u," and ",v,"\nare doubly-twisted
            conjugated by element ",w, "\n");
    fi;
    w := NextIterator(I);
  od;
  return fail;
end;
```

## A.2.2 "dtc.gap"

```
# File name: dtc.gap
# Author: Michele Feltz
# Last modified: 19/12/2009

#************************Terminology************************
# G,H: finitely generated free groups
# f,g: homomorphisms between finitely generated free groups
# u,v: freely reduced words
# l,m,n: positive integers
#***********************************************************

# definition of some free groups
F11 := FreeGroup("b1");
F12 := FreeGroup("b2");
F2 := FreeGroup("a", "b");
F3 := FreeGroup("a", "b", "c");
F4 := FreeGroup("a", "b", "c", "d");

# rel( G, l) generates a random element of the free group G of
# length at most l
rel := function( G, l)
  local letterlist, gens, el;
  gens := GeneratorsOfGroup(G);
  letterlist := function()
    return List([1..l], n -> Random(gens)^Random([-1,1]));
  end;
  el := Product(letterlist());
  return el;
end;

# FreeProductHomomorphism( f, u, v, G, H, GFF) computes the
# homomormorphism f*u*v : GFF -> H (where GFF is the free
# product G*F11*F12) defined as follows:
# a -> (af), b1 -> u, b2 -> v (with a generators of G,
# b1 and b2 generators of F11 and F12 respectively)
FreeProductHomomorphism := function( f, u, v, G, H, GFF)
    local a,l;
    l := Concatenation(List(GeneratorsOfGroup(G),a -> a^f),
```

```
            [u,v]);
      return GroupHomomorphismByImagesNC(GFF,H,
                                        GeneratorsOfGroup(GFF),
                                        l);
end;

# HasBSLW( f, g, u, v, G, H) indicates whether or not
# |Rem_{f*u*v}(a)| >= |(ag)| for all the generators a of G
HasBSLW := function( f, g, u, v, G, H)
  local a, r, l, p1, fp, fph;
  fp := FreeProduct(G,F11,F12);
  fph := FreeProductHomomorphism(f,u,v,G,H,fp);
  p1 := Embedding(fp,1);
  if not HasRemnant(fph,fp,H) then
  return false;
  else
  for a in GeneratorsOfGroup(G) do
      r := RemnantOfGenerator(fph,a^p1,fp,H);
      l := a^g;
      if (Length(r) < Length(l)) then
      return false;
      fi;
    od;
    return true;
  fi;
end;

# If HasBSL(f,g,G,H) returns true, then MinBSL( f, g, G, H)
# computes the minimum of the difference |Rem_{f}(a)| - |(ag)|
# over all the generators a of G, otherwise MinBSL returns fail
MinBSL := function( f, g, G, H)
  local a, r, l, l1, l2, d, m;
  if not HasBSL(f,g,G,H) then
    return fail;
  fi;
  m := 100;
  for a in GeneratorsOfGroup(G) do
    r := RemnantOfGenerator(f, a, G, H);
    l := a^g;
    l1 := Length(r);
    l2 := Length(l);
    d := l1-l2;
```

```
    if (d <= m) then
      m := d;
      fi;
    od;
  return m;
end;

# IsDTCBSL( f, g, u, v, G, H) provides an algorithm for
# deciding doubly-twisted conjugacy of the words u and v
# (in case HasBSL(f,g,G,H) returns true)
IsDTCBSL := function( f, g, u, v, G, H)
  local SB,l,b;
  if not HasBSL(f,g,G,H) then
    Print("BSL condition not fulfilled!\n");
  else
    if u = v then
      Print("The elements ", u, " and ", v, "\nare doubly-
             twisted conjugated by the identity element of ",
             G,"\n");
    elif HasBSLW(f,g,u,v,G,H) then
      Print("The Elements ", u," and ", v," are in different
             classes.\n");
    else
      l := MinBSL(f,g,G,H);
      SB := (Length(u) + Length(v))/l;
      b := IsDoublyTwistedConjugateBruteForce(f,g,u,v,G,H,SB);
      if b = true then
        return ;
      else
        return b;
      fi;
    fi;
  fi;
end;

# DTCBSLTests( f, g, G, H, m, n) generates random pairs of dis-
# tinct elements from the free group H of lengths at most m
# which satisfy the doubly-twisted conjugacy relation for some
# w in G (out of n trials) (in case HasBSL(f,g,G,H) returns
# true)
DTCBSLTests := function( f, g, G, H, m, n)
  local SB,l,b,i,u,v;
```

```
    if not HasBSL(f,g,G,H) then
      Print("BSL condition not fulfilled!\n");
      return ;
    fi;
    i := n;
    while i >= 0 do
      u := rel(G,m);
      v := rel(G,m);
      if not (u = v or HasBSLW(f,g,u,v,G,H)) then
        l := MinBSL(f,g,G,H);
        SB := (Length(u) + Length(v))/l;
        b := IsDoublyTwistedConjugateBruteForce(f,g,u,v,G,H,SB);
      fi;
      i := i-1;
    od;
end;

# definition of some free group homomorphisms and words
rh1 := GroupHomomorphismByImagesNC(F2,F2,GeneratorsOfGroup(F2),
       [F2.2^2*F2.1*F2.2^2*F2.1^-2,
        F2.1^-1*F2.2^-1*F2.1*F2.2^-1*F2.1^2*F2.2^-1]);
rh2 := GroupHomomorphismByImagesNC(F2,F2,GeneratorsOfGroup(F2),
       [F2.2^-1,F2.1*F2.2^-2]);

rh3 := GroupHomomorphismByImagesNC(F2,F2,GeneratorsOfGroup(F2),
       [F2.2^-3*F2.1*F2.2^-1,
        F2.1^-1*F2.2^-1*F2.1*F2.2*F2.1^-1]);
rh4 := GroupHomomorphismByImagesNC(F2,F2,GeneratorsOfGroup(F2),
       [F2.1,F2.2]);

u1 := F2.2*F2.1^-1*F2.2*F2.1^-1;
v1 := F2.2*F2.1^-1*F2.2^-2;

u2 := F2.1^-1*F2.2^-1;
v2 := F2.2*F2.1^-2*F2.2^-1;
```

## A.3 Tests

To give the reader an idea of how to use the previous GAP programs, we generate some tests on the command line. We make use of the homomorphisms and words that have been previously defined in the file "dtc.gap". Note that, since the second file depends on the first one, the latter needs to be loaded first.

```
gap> Read("GAPP3/final/rembsl.gap");
gap> Read("GAPP3/final/dtc.gap");
```

The next tests relate to the example in Section 3.11.

```
gap> HasRemnant(rh1,F2,F2);
true
gap> RemnantTriple(rh1,F2.1,F2,F2);
[ b, b*a*b^2*a^-2, <identity ...> ]
gap> RemnantTriple(rh1,F2.2,F2,F2);
[ <identity ...>, a^-1*b^-1*a*b^-1*a^2, b^-1 ]
gap> HasBSL(rh1,rh2,F2,F2);
true
gap> MinBSL(rh1,rh2,F2,F2);
3
gap> GFF:=FreeProduct(F2,F11,F12);;
gap> fph:=FreeProductHomomorphism(rh1,u1,v1,F2,F2,GFF);
[ f1, f2, f3, f4 ] -> [ b^2*a*b^2*a^-2, a^-1*b^-1*a*b^-1*a^2*b^-1,
  b*a^-1*b*a^-1, b*a^-1*b^-2 ]
gap> HasBSLW(rh1,rh2,u1,v1,F2,F2);
false
gap> IsDTCBSL(rh1,rh2,u1,v1,F2,F2);
fail
gap> HasBSLW(rh1,rh2,u2,v2,F2,F2);
false
gap> IsDTCBSL(rh1,rh2,u2,v2,F2,F2);
The elements a^-1*b^-1 and b*a^-2*b^-1
are doubly-twisted conjugated by element b
```

Now we generate random pairs of different elements from the free group $F_2$ that are doubly-twisted conjugated (out of 10000 trials).

```
gap> DTCBSLTests(rh1,rh2,F2,F2,6,10000);
The elements b^2*a*b and a^2*b^-2
are doubly-twisted conjugated by element a
The elements b^2 and a^2*b^-2*a^-1*b^-1
are doubly-twisted conjugated by element a
The elements a^-1*b^-1 and b*a^-2*b^-1
are doubly-twisted conjugated by element b
The elements a^-1*b^-1*a*b^-1 and b*a^-1*b^-2
are doubly-twisted conjugated by element b

gap> DTCBSLTests(rh1,rh2,F2,F2,8,10000);
The elements b^2*a*b^3 and a^2
are doubly-twisted conjugated by element a
```

Notice that $rh4$ is the identity homomorphism and that the hypotheses of Corollary 3.12 for $rh3, rh4$ are satisfied, so that the following tests generate pairs of different random elements which belong to the same *twisted* conjugacy class.

```
gap> HasBSL(rh3,rh4,F2,F2);
true
gap> DTCBSLTests(rh3,rh4,F2,F2,8,10000);
The elements a^-1*b^-1*a*b*a^-1*b and b^2
are doubly-twisted conjugated by element b
The elements b*a and b^-3*a
are doubly-twisted conjugated by element a^-1
The elements b*a and b^-3*a
are doubly-twisted conjugated by element a^-1
The elements a^-1*b^-1 and a*b^-1*a^-1*b
are doubly-twisted conjugated by element b
The elements a*b^-1*a*b and a^-1*b^-1*a^2
are doubly-twisted conjugated by element b^-1
The elements a^-1*b^-1*a*b and a*b
are doubly-twisted conjugated by element b
The elements b*a^-1*b*a^3 and b^-2*a^2
are doubly-twisted conjugated by element a^-1
The elements a*b^-1*a^-1*b*a^-1*b and a^-2
are doubly-twisted conjugated by element b^-1
The elements b^-3*a^-1*b^-1*a^-1 and b*a^-2*b^-1
are doubly-twisted conjugated by element a
The elements a*b^-1 and a*b^-1*a^-1*b*a^2
are doubly-twisted conjugated by element b

gap> DTCBSLTests(rh3,rh4,F2,F2,10,10000);
The elements b^-3*a^-1 and b*a^-1*b^4
are doubly-twisted conjugated by element a^2
The elements b^-1*a*b^-1*a and b^-3*a*b^-2*a*b^-1
are doubly-twisted conjugated by element a^-1
The elements b^-2*a*b^-1 and b*a^-1*b*a*b^-1*a
are doubly-twisted conjugated by element a
The elements b*a^-1*b*a^-1 and b^-2*a^-2
are doubly-twisted conjugated by element a^-1
The elements a*b^-1*a^-1*b*a^-1*b and a^-2
are doubly-twisted conjugated by element b^-1
The elements a^-1*b^-1*a*b^-1 and a*b^-1
are doubly-twisted conjugated by element b

gap> DTCBSLTests(rh3,rh4,F2,F2,12,10000);
```

```
The elements a^-1*b^-1*a*b*a^-1*b and b^2
are doubly-twisted conjugated by element b
The elements b*a^-1*b^4 and b*a^-1
are doubly-twisted conjugated by element a^-1
gap> quit;
```

# Appendix B

# The Protocol in GAP

In this appendix we give a user-friendly implementation of the Shpilrain-Ushakov authentication protocol described in Chapter 4. First we list the GAP code of the file "protocol.gap" in Subsection B.1, then we execute the file in Subsection B.2.

## B.1 Code listing

### B.1.1 "protocol.gap"

```
# File name: protocol.gap
# Authors: Sasa Radomirovic, Laura Ciobanu, Michele Feltz
# Last modified: 21/01/2010

# r is the ring of integers mod 2, that is r=GF(2)
r:=Integers mod 2;

# x is an indeterminate over r
x:=Indeterminate(r, "x");

# We truncate polynomials at degree N. If no value has been
# assigned to N yet, then we set N to 7
if First(NamesUserGVars(),t->t="N") = fail then
  N:=7;
fi;

# poly( vector) takes as input a vector [a_n, a_{n-1}, ..., a_0]
# of coefficients in r and returns the polynomial function
# x -> a_n x^n + a_{n-1} x^{n-1} + ... + a_0 x_0
poly:=function(vector)
  return u->Iterated(vector,function(a,b) return a*u+b; end);
```

```
end;

# poly( vector) generates a polynomial whose constant term is 0
# Note: the degree of the polynomial will be one higher than
# above for the same vector
hpoly:=function(vector)
  return u->u*Iterated(vector,function(a,b) return a*u+b; end);
end;

# randomVector( n, R) returns a random vector of length n of
# elements from R (where R can be r, GF(2), [0,1], ...)
randomVector:=function(n,R)
  local t;
  return List([1..n],t->Random(R));
end;

# randmat() returns a random matrix, whose entries
# are polynomial functions over r of degree N-1.
randmat:=function()
  local rv1,f1,rv2,f2,rv3,f3,rv4,f4;
  rv1:=randomVector(N,r);
  f1:=u->poly(rv1)(u);
  rv2:=randomVector(N,r);
  f2:=u->poly(rv2)(u);
  rv3:=randomVector(N,r);
  f3:=u->poly(rv3)(u);
  rv4:=randomVector(N,r);
  f4:=u->poly(rv4)(u);
  return([[f1,f2],[f3,f4]]);
end;

# mymod( a, b) computes a mod b
mymod:=function(a,b)
  if IsUnit(a) or IsUnit(a+Z(2)) then
    return a;
  else
    return a mod b;
  fi;
end;

# matrixMod( A, m) reduces all entries in the 2x2 matrix A by m
matrixMod:=function(A,m)
  return [[mymod(A[1][1],m),mymod(A[1][2],m)],[mymod(A[2][1],m),
          mymod(A[2][2],m)]];
```

```
end;

# polmat( mypol, mymat) apply polynomial mypol to matrix mymat
# whose entries are polynomial functions
polmat:=function(mypol, mymat)
 local m;
 m:=[[mymat[1][1](mypol),mymat[1][2](mypol)],[mymat[2][1](mypol),
       mymat[2][2](mypol)]];
 return matrixMod(m,x^N);
end;

# ldeg( pol) computes the degree of the lowest degree monomial
# in pol; returns N+1 as the "degree" of the zero polynomial
ldeg:=function(pol)
  local i;
  i:=1;
  if pol = (x-x) then return N+1; fi;
  if pol = 0*Z(2) then return N+1; fi;
  while(mymod(pol,x^i)=(x-x)) do i:=i+1; od;
  return i-1;
end;

# polyvec( p) takes as input a polynomial p and returns
# a vector v such that poly(v)(x)=p
polyvec:=function(p)
  local v,q,pld;
  if ldeg(p)>N then return([0]); fi;
  v:=[]; q:=p; pld:=0;
  if ldeg(p)>0 then v:=[0]; fi;
  while ldeg(q)<N do
    v:=Concatenation([1],Concatenation(List([1..ldeg(q)-pld-1],
                                            t->0),v));
    pld:=ldeg(q);
    q:=q-x^pld;
  od;
  return(v);
end;

# listmat( mat) turns a regular matrix (whose entries are
# polynomials) into a matrix whose entries are vectors
listmat:=function(mat)
  return([polyvec(mat[1][1]),polyvec(mat[1][2]),
          polyvec(mat[2][1]),polyvec(mat[2][2])]);
end;
```

```
# funcmat( ma) turns a regular matrix into a matrix whose
# entries are functions
funcmat:=function(ma)
  local mat;
  mat:=listmat(ma);
  return([[poly(mat[1]),poly(mat[2])],[poly(mat[3]),
          poly(mat[4])]]);
end;

#********************* SECRET KEY ***************************
# the secret key is represented as a matrix of 4 random
# polynomial functions (the polynomials are of degree
# at most N-1 with constant term 1 (that is, Z(2))
rv1:=randomVector(N-1,r);
s1:=u->hpoly(rv1)(u)+Z(2);

rv2:=randomVector(N-1,r);
s2:=u->hpoly(rv2)(u)+Z(2);

rv3:=randomVector(N-1,r);
s3:=u->hpoly(rv3)(u)+Z(2);

rv4:=randomVector(N-1,r);
s4:=u->hpoly(rv4)(u)+Z(2);

# the secret!
secret:=[[s1,s2],[s3,s4]];

# reveals allows to see the secret, more precisely, the
# polynomial entries (not the corresponding polynomial
# functions)
reveals:=polmat(x,secret);

#********************* PUBLIC KEY ***************************
# let us generate random endomorphisms of degree at most N-1
# over r[x]; the endomorphisms are represented as polynomials
# with 0 constant term
rphi:=randomVector(N-1,r);
phi:=hpoly(rphi)(x);

rpsi:=randomVector(N-1,r);
```

```
psi:=hpoly(rpsi)(x);

# another public key w is represented as a matrix whose entries
# are random polynomials of degree at most N-1 with constant
# term 1 (that is, Z(2))
rv5:=randomVector(N-1,r);
w1:=u->hpoly(rv5)(u)+Z(2);

rv6:=randomVector(N-1,r);
w2:=u->hpoly(rv6)(u)+Z(2);

rv7:=randomVector(N-1,r);
w3:=u->hpoly(rv7)(u)+Z(2);

rv8:=randomVector(N-1,r);
w4:=u->hpoly(rv8)(u)+Z(2);

w:=polmat(x,[[w1,w2],[w3,w4]]);

#----------------- MOST IMPORTANT FUNCTION --------------------
# doubleTwist( s, mat) computes the doubly-twisted conjugacy
# relation on mat, that is, (Transpose((sphi)))*mat*(spsi)
doubleTwist:=function(s,mat)
  return(matrixMod(TransposedMat(polmat(phi,s))*mat*
                                polmat(psi,s),x^N));
end;
#--------------------------------------------------------------

# the last public key t is obtained by doubly-twisting w
t:=doubleTwist(secret,w);

#*************** START PROTOCOL: SINGLE ROUND *****************
# Alice's choice of random matrix r=rm:
rm:=randmat();
revealr:=polmat(x,rm);

# Alice computes sr from r=rm and her secret key s=secret
revealsr:=matrixMod(polmat(x,secret)*polmat(x,rm),x^N);
sr:=funcmat(revealsr);

# Visualize parameters, public key and product sr=secret*rm
Print("Polynomials are truncated at N = ",N,"\n");
Print("Public Key is\nphi = ",phi,"\npsi = ",psi,"\n  w = ",w,"\n
```

```
        t = ",t,"\n");
Print("Product of secret with random matrix is = ",revealsr,"\n");
# Print("secret key: ",reveals,"\n");

Print("Alice choses at random r = ",revealr,"\n");

u:=doubleTwist(rm,t);
Print("Alice sends commitment u:=doubleTwist(r,t) = ",u,"\n");
Print("Bob chooses a random bit c and sends it as a challenge
        to Alice.\n");
v:=rm;
Print("Suppose c=0, then Alice sends to Bob v=r = ",revealr,"\n");
Print("Bob checks whether or not u=doubleTwist(v,t): ",
        u=doubleTwist(v,t),"\n");

v:=sr;
Print("Suppose c=1, then Alice sends to Bob v=sr = ",revealsr,"\n");
Print("Bob checks whether or not u=doubleTwist(v,w): ",
        u=doubleTwist(v,w),"\n");
```

## B.2 Test

When executing the file "protocol.gap" we obtained the following output. Note that reloading the file yields different public and secret keys, thus a different protocol run.

```
gap> Read("GAPP3/Crypto/protocol.gap");
Polynomials are truncated at N = 7
Public Key is
phi = x^3+x
psi = x^6+x^2+x
  w = [ [ x^4+Z(2)^0, x^4+x+Z(2)^0 ],
  [ x^6+x^3+x+Z(2)^0, x^4+x^3+x^2+x+Z(2)^0 ] ]
 t = [ [ x^5+x^4+x^2+x, x^6+x ], [ x^6+x^5+x^3+x, x^5+x ] ]
Product of secret with random matrix is =
[ [ x^6+x^3+Z(2)^0, x^2+x+Z(2)^0 ],
  [ x^6+x^5+x^3+x^2+Z(2)^0, x^6+x^5+x^4+x^3+Z(2)^0 ] ]
Alice choses at random r = [ [ x^6+x^5+x^4+Z(2)^0, x^3+x ],
  [ x^6+x^3+x, x^5+Z(2)^0 ] ]
Alice sends commitment u:=doubleTwist(r,t) =
[ [ x^5+x^2+x, x^4+x^3+x ], [ x^6+x^5+x, x^6+x^5+x ] ]
Bob chooses a random bit c and sends it as a challenge to Alice.
Suppose c=0, then Alice sends to Bob v=r =
[ [ x^6+x^5+x^4+Z(2)^0, x^3+x ], [ x^6+x^3+x, x^5+Z(2)^0 ] ]
```

```
Bob checks whether or not u=doubleTwist(v,t): true
Suppose c=1, then Alice sends to Bob v=sr =
[ [ x^6+x^3+Z(2)^0, x^2+x+Z(2)^0 ],
  [ x^6+x^5+x^3+x^2+Z(2)^0, x^6+x^5+x^4+x^3+Z(2)^0 ] ]
Bob checks whether or not u=doubleTwist(v,w): true
```